\documentclass[a4paper]{article}

\usepackage[english]{babel}

\usepackage[utf8x]{inputenc}
\usepackage[T1]{fontenc}

\usepackage{multirow, longtable, caption, booktabs}
\usepackage{pdflscape}

\usepackage[a4paper,top=3cm,bottom=2cm,left=3cm,right=3cm,marginparwidth=1.75cm]{geometry}

\usepackage{amsmath}
\usepackage{mathrsfs}
\usepackage{amssymb}
\usepackage{graphicx}
\usepackage{subfigure}
\usepackage[colorinlistoftodos]{todonotes}
\usepackage[colorlinks=true, allcolors=blue]{hyperref}
\usepackage{setspace}

\newcommand{\cdummy}{\cdot}

\usepackage{dsfont, makecell}
\usepackage[authoryear, round]{natbib}

\usepackage{amsthm}
\newtheorem{theorem}{Theorem}
\newtheorem{proposition}{Proposition}

\DeclareMathOperator*{\argmin}{argmin}

\global\long\def\prox{\operatornamewithlimits{prox}}%
\usepackage[ruled, linesnumbered]{algorithm2e}

\usepackage{bm}

\newcommand{\bx}{\mathbf{x}}
\newcommand{\bb}{\mathbf{b}}
\newcommand{\bc}{\mathbf{c}}
\newcommand{\by}{\mathbf{y}}

\newcommand{\bs}{\mathbf{s}}
\newcommand{\bq}{\mathbf{q}}
\newcommand{\bp}{\mathbf{p}}
\newcommand{\br}{\mathbf{r}}
\newcommand{\bu}{\mathbf{u}}
\newcommand{\bv}{\mathbf{v}}
\newcommand{\bA}{\mathbf{A}}
\newcommand{\bD}{\mathbf{D}}
\newcommand{\bQ}{\mathbf{Q}}
\newcommand{\bI}{\mathbf{I}}
\newcommand{\bzero}{\mathbf{0}}
\newcommand{\bz}{\mathbf{z}}
\newcommand{\bP}{\mathbf{P}}
\newcommand{\bLambda}{\bm{\Lambda}}
\newcommand{\bX}{\mathbf{X}}
\newcommand{\bE}{\mathbf{E}}
\newcommand{\bM}{\mathbf{M}}
\newcommand{\bh}{\mathbf{h}}
\newcommand{\bL}{\mathbf{L}}
\newcommand{\bw}{\mathbf{w}}
\newcommand{\bt}{\mathbf{t}}
\newcommand{\bK}{\mathbf{K}}
\newcommand{\bG}{\mathbf{G}}
\newcommand{\bH}{\mathbf{H}}
\newcommand{\bF}{\mathbf{F}}
\newcommand{\bzeta}{\bm{\zeta}}
\newcommand{\bxi}{\bm{\xi}}

\title{An Enhanced ADMM-based Interior Point Method for Linear and Conic Optimization}

\author{
Qi Deng$^{1}$\and 
Qing Feng$^{2}$ \and 
Wenzhi Gao$^{3}$\and
Dongdong Ge$^{1}$\and
Bo Jiang$^{4}$\footnote{Correspondence to \url{jiang.bo@mail.shufe.edu.cn}.}\and
Yuntian Jiang$^{4}$\and
Jingsong Liu$^{4}$\and
Tianhao Liu$^{4}$\and
Chenyu Xue$^{4}$\and
Yinyu Ye$^{3}$\and
Chuwen Zhang$^{4}$\and \\
$^{1}$Shanghai JiaoTong University, $^{2}$Cornell University, \\
$^{3}$Stanford University, $^{4}$Shanghai University of Finance and Economics
}
\begin{document}
\maketitle

\begin{abstract}
    The ADMM-based interior point (ABIP, Lin et al.\ 2021) method is a hybrid algorithm that effectively combines interior point method (IPM) and first-order methods to achieve a performance boost in large-scale linear optimization. Different from traditional IPM that relies on computationally intensive Newton steps, the ABIP method applies the alternating direction method of multipliers (ADMM) to approximately solve the barrier penalized problem. However, similar to other first-order methods, this technique remains sensitive to condition number and inverse precision.
In this paper, we provide an enhanced ABIP method with multiple improvements. Firstly, we develop an ABIP method to solve the general linear conic optimization and establish the associated iteration complexity.  
Secondly, inspired by some existing methods, we develop different implementation strategies for ABIP method, which substantially improve its performance in linear optimization.
Finally, we conduct extensive numerical experiments in both synthetic and real-world datasets to demonstrate the empirical advantage of our developments. In particular, the enhanced ABIP method achieves a 5.8x reduction in the geometric mean of run time on $105$ selected LP instances from Netlib, and
it exhibits advantages in certain structured problems such as 
SVM and PageRank.
However, the enhanced ABIP method still falls behind commercial solvers in many benchmarks, especially when high accuracy is desired. We posit that it can serve as a complementary tool alongside well-established solvers.
\end{abstract}

\noindent {\bf Keywords:}
Linear optimization; Conic optimization; ADMM; Interior point method; Implementation improvement; Iteration complexity

\section{Introduction}
In this paper, we consider the following linear conic optimization problem with standard primal (P) and dual (D) forms:
\begin{equation}
\label{prob:conic-main}
    \text{(P)}
	\begin{aligned}
	   \min & \quad \mathbf{c}^{T}\mathbf{x} \\
	   \quad	\text{s.t.} &\quad \mathbf{A} \mathbf{x} = \mathbf{b} \\
	   &\quad \bx \in \mathcal{K}
	\end{aligned}
    \quad\quad
    \text{(D)}
    \begin{aligned}
	   \max & \quad \mathbf{b}^{T}\mathbf{y} \\
	   \quad	\text{s.t.} &\quad \mathbf{A}^{T} \mathbf{y} + \mathbf{s} = \mathbf{c} \\
	   &\quad \mathbf{s} \in \mathcal{K}^{*}
	\end{aligned}
\end{equation}   
where $\bc\in\mathbb{R}^{n}$, $\bb\in\mathbb{R}^{m}$, $\bA\in\mathbb{R}^{m\times n}$, $\mathcal{K}$ is a nonempty, closed, and convex cone with its dual cone defined as $\mathcal{K}^{*} = \left\lbrace \by: \bx^{T}\by \geq 0,\forall \bx\in\mathcal{K} \right\rbrace$. In fact, this problem has strong modeling power, which includes linear optimization (LP), second-order cone optimization (SOCP), and semidefinite optimization (SDP) as special cases. 
Moreover, many important problems, such as quadratic optimization (QP) and quadratically constrained QP (QCQP), can be translated to equivalent conic formulations described by \eqref{prob:conic-main} \citep{alizadeh2003second}.

Due to the wide applications in engineering and data science, developing efficient and accurate algorithms for linear and conic optimization has been a central topic in the optimization field in the past few decades. The traditional method for solving problem \eqref{prob:conic-main} is interior point method (IPM), which resorts to
a sequence of log-barrier penalty subproblems and requires one step of Newton's method to solve each subproblem \citep{nesterov1994interior}. Open-source and commercial solvers based on IPM, such as SeDuMi~\citep{sturm1999using}, SDPT3~\citep{toh1999sdpt3}, MOSEK~\citep{aps2019mosek}, GUROBI~\citep{gurobi}, COPT~\citep{ge2022cardinal}, are well-developed and have received great success in practice. Recently, new variants of IPM have been proposed \citep{pougkakiotis2022interior, cipolla2023proximal}. These algorithms show promising improvements compared to the traditional IPM, and their algorithm frameworks still fall in the scope of solving regularized Newton equations.
However, despite the fact that IPM can achieve fast convergence to high-accuracy solution, the computational cost can be a major concern, as solving a sequence of Newton equations can be highly expensive for large-scale or high-dimensional problems.

Compared to IPM, first-order methods are considered to be more scalable due to the low per-iteration cost and avoidance of solving Newton's equation. Initial attempts in this direction start by replacing the matrix decomposition with iterative methods, and it is shown that convergence of IPM can be preserved even if the Newton system is not accurately solved \citep{zhou2004polynomiality, bellavia2004convergence, lu2006iterative, al2009convergence, zanetti2023new}.
Lately, there has been a growing interest in developing first-order methods for solving large-scale LP or conic optimization \citep{YangST15, SCS-2016, PDHG-2021}. In particular, \cite{SCS-2016} develop the Splitting Conic Solver (SCS) for general conic LP, which applies the alternating direction method of multipliers (ADMM) \citep{BoydPCPE11} to solve the homogeneous self-dual (HSD) reformulation \citep{YeToddMizuno-1994} of the conic problem \eqref{prob:conic-main}. 
The numerical results in \cite{SCS-2016} showcase the superior performance of SCS over traditional IPM for several large-scale conic problems.  \cite{sopasakis2019superscs} present a new Douglas-Rachford splitting method for solving the HSD system, which uses the quasi-Newton directions, such as restarted Broyden directions, and Anderson’s acceleration to further improve the convergence performance. 

For LP problems, \cite{ABIP-2021} propose the ADMM-based Interior Point (ABIP) method, which can be viewed as a hybrid algorithm of the path-following IPM and ADMM. More specifically, it constructs a sequence of HSD reformulation of the LP problem with diminishing log-barrier parameters and uses ADMM to approximately solve each subproblem associated with a fixed log-barrier parameter. Therefore, it is expected to inherit some merits from both methods. For ease of exposition, we use ``ABIP'' and ``the ABIP method'' interchangeably in the rest of the paper.
Very recently, \cite{PDHG-2021} apply the primal-dual hybrid gradient (PDHG) \citep{ChambolleP11} method to solve LP based on its saddle point formulation. \cite{PDLP-2021} further propose a practical first-order method for LP (PDLP), which is an enhanced version of PDHG by combining several advanced implementation techniques. It is shown in \cite{PDLP-2021} that PDLP even outperforms a commercial LP solver in a large-scale application, i.e., the PageRank problem. 

While earlier papers consider linear optimization, some recent effort has been made in developing first-order methods for QP.
\cite{stellato2020osqp} present OSQP, a convex quadratic optimization solver based on ADMM, which is competitive against commercial solvers. \cite{garstka2021cosmo} propose a new conic operator splitting method (COSMO) solver, which extends ADMM to deal with more general conic-constrained QP. By using chordal decomposition and some new clique merging techniques, they significantly improve the algorithm performance on solving some large-scale SDP. \cite{o2021operator} extend the SCS framework~\citep{SCS-2016} to solve quadratic conic optimization based on a more general formulation of the linear complementarity problem. The new implementation \citep{scs} has shown a solid empirical advantage in the infeasible problem while preserving great efficiency in feasible problems.

\paragraph{Contributions.}
In this paper, we continue the development of ABIP in \cite{ABIP-2021} along several new directions. Our contributions can be summarized as follows. 

First, we present an important extension of ABIP such that the new solver can directly handle more general conic constraints. Theoretically, we show that the extended ABIP obtains an $\mathcal{O}\left(\left({1}/{\epsilon}\right)\log\left(1/\epsilon\right)\right)$ complexity bound, thus generalizing the complexity results previously known for LP \citep{ABIP-2021} to a more general conic setting. For practical implementation, we show that the proximal operator associated with the log-barrier subproblem in ABIP can be efficiently computed. 
For several important applications in machine learning, such as LASSO and SVM, we develop customized linear system solvers to further accelerate ABIP for specific large-scale problems. 
We use extensive experiments on both synthetic and real-world datasets to show that enhanced ABIP compares favorably against many popular open-source and commercial solvers. 

Second, we significantly improve the practical performance of ABIP for LP by developing several acceleration strategies. Those strategies were mostly motivated and adapted from the techniques used in the previous literature, including the adaptive strategy for choosing barrier parameter $\mu$, restart scheme, new inner loop convergence criteria, half update in ADMM, Presolve and preconditioning, and a tailored acceleration strategy for null objective problems. 
We further propose a decision tree-based approach to efficiently integrate those new strategies. 
With all the mentioned acceleration techniques, our numerical experiments show that the enhanced ABIP achieves a 5.8x reduction in the geometric mean of run time on $105$ selected problem instances from Netlib. 
In large-scale applications such as the (staircase) PageRank problem, the enhanced ABIP is highly competitive and often outperforms PDLP \citep{PDLP-2021}.
It is worth noting that the enhanced ABIP still falls behind commercial solvers across a variety of benchmarks, especially when a highly accurate solution is desired. This limitation seems inherent to all first-order methods, and addressing it remains a subject for future research.

\paragraph{Organizing the paper.}
This paper proceeds as follows. Section~\ref{Sec:conic} generalizes ABIP to solve the convex conic optimization.  Section~\ref{Sec:new-strategy} develops new strategies to further accelerate the practical performance of ABIP for LP. Section~\ref{Sec:numerical} conducts a detailed and extensive experimental study to demonstrate the empirical advantage of ABIP. 

\paragraph{Notation and terminology.}
We use bold-face letters to denote matrices (i.e., $\bA, \bQ$) and vectors (i.e., $\bx, \by$). Let $\bzero$ be a  vector of zeros, with its dimensionality unspecified, whenever it is clear from the context. Let $\mathbb{R}^n$ be the $n$-dimensional Euclidean space.
We use $\bx \ge \by$ to denote the element-wise inequality $x_i\ge y_i$ for all $i$. The nonnegative orthant is defined by $\mathbb{R}^n_+=\{ \bx \mid \bx \ge \bzero, \bx \in \mathbb{R}^n \}$. The regular second-order cone and the rotated second-order cone are denoted by $\text{SOC}_{1+n}=\{(t, \by)\in\mathbb{R}_+\times\mathbb{R}^{n}: \|\by\| \leq t\}$ and $\text{RSOC}_{2+n} = \{(t,s,\by)\in\mathbb{R}^{2}_+\times\mathbb{R}^{n}:ts\geq \frac{1}{2}\|\by\|_2^2\}$. respectively.

\section{ABIP for conic optimization}
\label{Sec:conic}
The original ABIP is designed for LP only. In this section, we extend it to solve general linear conic optimization problems where $\mathcal{K}$ can be represented by the Cartesian product of different general convex cones. 

\paragraph{Motivating problems.} We consider the following convex quadratic problem: 
\begin{equation}\label{prob:qp}
    \begin{aligned}
    \min_{\bz} &\quad \frac{1}{2} \bz^T \bP\bz + \bq^T \bz \\
    \text{s.t.} & \quad \bar{\bA} \bz = \bar{\bb} \\
                & \quad \bz \geq \bzero
    \end{aligned}  
\end{equation}
where $\bar{\bA} \in \mathbb{R}^{m\times n}$, and $\bP\in\mathbb{R}^{n\times n}$ is a positive semi-definite matrix. 
It has extensive applications in machine learning and data science, such as SVM \citep{cortes1995support}, LASSO \citep{tibshirani1996regression}, and portfolio optimization \citep{markowitz1991foundations}. Due to the existence of the quadratic term in the objective function, problem~\eqref{prob:qp} can not be directly solved by ABIP. However, it can be reformulated into a linear conic optimization problem. 

Let $\bP = \bLambda^{T}\bLambda$ be the factorization where $\bLambda\in\mathbb{R}^{r\times n}$, 
problem~\eqref{prob:qp} can be expressed as the following problem:
\begin{equation*}
	\begin{aligned}
		\min &\quad \frac{1}{2}\bar{\bx}^{T}\bar{\bx} + \bq^{T}\bar{\bz} \\
		\text{s.t.} &\quad\bar{\bA}\bar{\bz} = \bar{\bb} \\
		&\quad \bar{\bx} = \Lambda \bar{\bz} \\
		&\quad \bar{\bz}\geq \bzero
	\end{aligned}
\end{equation*}
which can be further transformed into the conic form:
\begin{equation*}
    \begin{aligned}
    	\min &\quad \nu + \bq^{T}\bar{\bz} \\
    	\text{s.t.} & \quad\eta = 1 \\
    	& \quad\bar{\bA}\bar{\bz} = \bar{\bb} \\
    	& \quad\bar{\bx} = \bLambda \bar{\bz} \\
    	& \quad\eta\nu \geq \frac{1}{2}\bar{\bx}^{T}\bar{\bx} \\
    	& \quad\bar{\bz} \geq \bzero
    \end{aligned}
    \ \Leftrightarrow\ 
    \begin{aligned}
        \min &\quad 
        \begin{bmatrix}
            0 \\ 1 \\ \bzero_{r\times 1} \\ \bq
        \end{bmatrix}^{T}
        \begin{bmatrix}
            \eta \\ \nu \\ \bar{\bx} \\ \bar{\bz}
        \end{bmatrix} \\
        \text{s.t.} & \quad
        \begin{bmatrix}
            1 \\
              & \bzero_{r\times 1} & \bI_{r} & -\bLambda \\
              &   &   & \bar{\bA}
        \end{bmatrix}
        \begin{bmatrix}
            \eta \\ \nu \\ \bar{\bx} \\ \bar{\bz}
        \end{bmatrix} = 
        \begin{bmatrix}
            1 \\ \bzero_{r\times 1} \\ \bar{\bb}
        \end{bmatrix} \\
        &\quad (\eta, \nu, \bar{\bx}) \in \text{RSOC}_{2+r}, \bar{\bz} \in \mathbb{R}_{+}^{n}
    \end{aligned}
\end{equation*}
where $\bI_{r}$ is $r$-dimensional identity matrix.

It should be noted that, although we use convex quadratic optimization as a motivating example, the conic formulation~\eqref{prob:conic-main}, as well as our subsequent algorithmic development and convergence analysis, still hold for general convex conic optimization such as SDP. We further provide an extensive numerical study under the SOCP setting in Section \ref{Sec:numerical}, while a more sophisticated empirical discussion beyond this setting will be left for future work.

\paragraph{HSD embedding for conic optimization.}
Before extending ABIP to linear conic optimization, we briefly describe the HSD embedding technique for convex conic optimization~\citep{Luo-1997, Zhang2004}. Given an initial solution  $\bx^{0}\in\text{int}\:\mathcal{K}$, $\bs^{0}\in\text{int}\:\mathcal{K}^{*}$ and $\by^{0}\in\mathbb{R}^{m}$, we can derive the following HSD embedding formulation:
\begin{equation}
	\begin{aligned}
		\min &\quad\left((\bx^{0})^{T}\bs^{0} + 1 \right) \theta \\
		\text{s.t.} &\quad \bA\bx -\bb\tau + \br_{p}\theta = \bzero \\
		&\quad -\bA^{T}\by + \bc\tau + \br_{d} \theta = \bs \\
		&\quad \bb^{T}\by - \bc^{T}\bx + r_{g}\theta = \kappa \\
		&\quad -\br_{p}^{T}\by - \br_{d}^{T}\bx - r_{g}\tau = -(\bx^{0})^{T}\bs^{0} - 1 \\
		&\quad \bx\in\mathcal{K}, \bs\in\mathcal{K}^{*}, \tau,\kappa\geq 0
	\end{aligned} \label{hsd} %
\end{equation}	
where $\br_{p} = \bb - \bA\bx^{0}$, $\br_{d} = \bs^{0} - \bc + \bA^{T}\by^{0}$, $r_{g} = 1 + \bc^{T}\bx^{0} - \bb^{T}\by^{0}$.
It is easy to check that the constraint matrix of~\eqref{hsd} is skew-symmetric, and hence the above problem is self-dual. 
The conic HSD embedding formulation~\eqref{hsd} enjoys several attractive properties~\citep{Zhang2004, Luo-1997}, which we list without proof as follows.
\begin{proposition}
\label{prop:feasible}
Problem \eqref{hsd} is strictly feasible. Specifically, there is a strictly feasible solution $(\by,\bx,\tau,\theta,\bs,\kappa)$  such that $\by = \by^{0}$, $\bx=\bx^{0}$, $\tau=1$, $\theta=1$, $\bs=\bs^{0}$, and $\kappa=1$.	
\end{proposition}

\begin{proposition}
\label{prop:infeasible}
Problem \eqref{hsd} has a maximally complementary optimal solution, denoted by $\left( \by^{*}, \bx^{*}, \tau^{*}, \theta^{*}, \bs^{*}, \kappa^{*}\right) $, such that $\theta^{*}=0$ and $(\bx^{*})^{T}\bs^{*}+\tau^{*}\kappa^{*} = 0$.
Moreover, the following statements hold: 
\begin{itemize}
    \item If $\tau^{*}>0$, then $\bx^{*} / \tau^{*} $ is an optimal solution for~\textup{(P)}, and $( \by^{*} / \tau^{*} , \bs^{*} / \tau^{*}) $ is an optimal solution for  \textup{(D)}.
    \item If $\kappa^{*}>0$, then either $\bc^{T}\bx^{*}<0$ or $\bb^{T}\by^{*}>0$. In the former case \textup{(D)} is infeasible, and in the latter case~\textup{(P)} is infeasible.
\end{itemize}
\end{proposition}

The first proposition shows that one can construct a strictly feasible solution to problem \eqref{hsd} based on a feasible solution to \eqref{prob:conic-main}. The second proposition shows that, one can obtain the optimal solution to \eqref{prob:conic-main} by solving \eqref{hsd}. Meanwhile, the HSD embedding technique can detect the infeasibility of a convex conic optimization problem, which is also inherited by ABIP. We remark that there are some other first-order algorithms that can handle infeasibility \citep{banjac2019infeasibility, applegate2021infeasibility}, which go beyond our scope in this paper since we mainly focus on how to solve the feasible problem efficiently.

\subsection{Algorithmic development}
In this subsection, we present the algorithm designed for linear conic optimization. Inspired by ABIP, we first rewrite problem \eqref{hsd} as:
\begin{equation}\label{hsd2}
	\begin{aligned}
		\min&\quad \beta((\bx^{0})^{T}\bs^{0}+1)\theta + \mathds{1}(\br=\bzero) + \mathds{1}(\xi = - (\bx^{0})^{T}\bs^{0} - 1) \\
		\text{s.t.}&\quad \bQ\bu = \bv \\
		&\quad (\bx,\bs,\tau,\kappa) \in \mathcal{K}\times \mathcal{K}^*\times\mathbb{R}_+\times \mathbb{R}_+, \ \by, \theta\ \text{free}, 
	\end{aligned}
\end{equation}
where $\bQ$, $\bu$, and $\bv$ are defined by, 
\[
\bQ = \begin{bmatrix}
	\bzero_{m\times m} & \bA & -\bb & \br_{p} \\
	-\bA^{T} & \bzero_{n\times n} & \bc & \br_{d} \\
	\bb^{T} & -\bc^{T} & 0 & \br_{g} \\
	-\br_{p}^{T} & -\br_{d}^{T} & -\br_{g} & 0
\end{bmatrix}, \bu = \begin{bmatrix} \by \\ \bx \\ \tau \\ \theta\end{bmatrix}, \,\,\text{and}\,\, \bv = \begin{bmatrix}\br \\ \bs \\ \kappa \\ \xi\end{bmatrix}.
\]
Then, we follow the framework of IPM and introduce the log-barrier penalty with parameter $\mu>0$ to the conic constraints, and move it to the objective function. Specifically, let $F(\bx)$ and $G(\bs)$ be the log-barrier functions associated with $\bx \in \mathcal{K}$ and $\bs \in \mathcal{K}^{*}$, respectively, 
and define
\begin{equation}\label{eq:B-conic}
    \begin{aligned}
        B(\bu,\bv,\mu) & = \beta((\bx^{0})^{T}\bs^{0}+1)\theta + \mathds{1}(\br=\bzero) + \mathds{1}(\xi = - (\bx^{0})^{T}\bs^{0} - 1) \\
        &\quad + \mu F(\bx) + \mu G(\bs) - \mu\log\tau - \mu\log\kappa.
    \end{aligned}
\end{equation}

Following \cite{ABIP-2021}, our enhanced ABIP has a double-loop structure. The details of the outer loop and the inner loop are discussed as follows. 

\paragraph{Outer loop.} At each outer iteration $k$, ABIP constructs the following log-barrier penalized problem with respect to $\mu^{(k)}$:
\begin{equation}
\begin{aligned}
    \min\ & B\left(\bu,\bv,\mu^{(k)}\right) \\
    \text{s.t.}\ & \bQ\bu=\bv.
\end{aligned} \label{prob:conic-barrier}
\end{equation}
Here, the sequence of $\{\mu^{(k)}\}$ is a diminishing sequence, i.e., $\lim_{k\rightarrow \infty}\mu^{(k)}=0$. In theoretical analysis, we let $\mu^{(k+1)} = \phi \mu^{(k)}$ with $\phi \in (0, 1)$.

\paragraph{Inner loop.} In the inner loop associated with $\mu^{(k)}$, ABIP use ADMM to approximately solves problem \eqref{prob:conic-barrier}. With variable splitting, \eqref{prob:conic-barrier} is equivalent to the following problem:
\begin{equation}\label{conic-subproblem}
	\begin{aligned}
		\min\ & \mathds{1}(\bQ\tilde{\bu}=\tilde{\bv}) + B\left(\bu,\bv,\mu^{(k)}\right) \\
		\text{s.t.}\ & (\tilde{\bu},\tilde{\bv}) = (\bu,\bv).
	\end{aligned} 
\end{equation}
Denote the augmented Lagrangian function by
\begin{align*}	\mathcal{L}_\beta(\tilde{\bu},\tilde{\bv},\bu,\bv,\mu^{(k)},\bp,\bq) = & \mathds{1}(\bQ\tilde{\bu} = \tilde{\bv}) + B\left(\bu,\bv,\mu^{(k)}\right) \\
& - \langle \beta(\bp,\bq), (\tilde{\bu},\tilde{\bv}) - (\bu,\bv)\rangle + \frac{\beta}{2}\|(\tilde{\bu},\tilde{\bv}) - (\bu,\bv)\|^{2},
\end{align*}
where $\beta>0$ is the same parameter as the one in~\eqref{hsd2}, and $\bp$, $\bq$ are the Lagrangian multipliers associated with the linear constraints. Then, in the $i$-th iteration, the update rule of ADMM is as follows:
\begin{align}
	(\tilde{\bu}^{(k)}_{i+1},\tilde{\bv}^{(k)}_{i+1}) &= \argmin_{\tilde{\bu},\tilde{\bv}} \mathcal{L}_\beta(\tilde{\bu},\tilde{\bv},\bu^{(k)}_{i},\bv^{(k)}_{i},\mu^{(k)},\bp^{(k)}_{i},\bq^{(k)}_{i}) \notag \\
 & = {\large \Pi}_{\bQ\bu=\bv}(\bu^{(k)}_{i}+\bp^{(k)}_{i}, \bv^{(k)}_{i}+\bq^{(k)}_{i}), \label{admm1} \\
	(\bu^{(k)}_{i+1}, \bv^{(k)}_{i+1}) &= \argmin_{\bu,\bv} \mathcal{L}_\beta(\tilde{\bu}^{(k)}_{i+1},\tilde{\bv}^{(k)}_{i+1},\bu,\bv,\mu^{(k)},\bp^{(k)}_{i},\bq^{(k)}_{i}), \label{admm2} \\
	(\bp^{(k)}_{i+1}, \bq^{(k)}_{i+1}) &= (\bp^{(k)}_{i}, \bq^{(k)}_{i}) - (\tilde{\bu}^{(k)}_{i+1}, \tilde{\bv}^{(k)}_{i+1}) + (\bu^{(k)}_{i+1}, \bv^{(k)}_{i+1}), \label{admm3}
\end{align}
where $\Pi_{\mathcal{S}}(\bx)$ denotes the Euclidean projection of $\bx$ onto the set $\mathcal{S}$. The above update rule seems to be complicated at the first glance. However, \cite{ABIP-2021} show that, for LP, the ADMM iteration can be substantially simplified by a customized initialization technique. We prove in the following theorem that this attractive feature still holds for linear conic optimization. The proof can be easily extended from \cite{ABIP-2021}, and we provide it in Appendix \ref{proof-thm-pq} for self-completeness.
\begin{theorem} 
\label{thm-pq}
At the $k$-th outer iteration of \textup{ABIP}, suppose that we initialize $\bp^{(k)}_{0} = \bv^{(k)}_{0}$ and $\bq^{(k)}_{0} = \bu^{(k)}_{0}$. Then, it holds that $\bp^{(k)}_{i}=\bv^{(k)}_{i}$ and $\bq^{(k)}_{i}=\bu^{(k)}_{i}$, for all $i\geq 0$.
\end{theorem}
Theorem~\ref{thm-pq} implies that, by carefully choosing the initial point at the beginning of the inner loop, we can eliminate the dual variables $\bp^{(k)}, \bq^{(k)}$ safely. Moreover, following a similar argument of \cite{ABIP-2021}, one can use the skew-symmetry property of $\bQ$ to obtain a simpler update for $\tilde{\bu}^{(k+1)}$, and use \eqref{admm3} to eliminate the update of $\tilde{\bv}_{i+1}^{(k)}$. 
In a nutshell, the above update rule can be simplified as:
\begin{align}
	\tilde{\bu}^{(k)}_{i+1} &= (\bI+\bQ)^{-1}(\bu^{(k)}_{i}+\bv^{(k)}_{i}), \label{cp-proj}\\
	\bu^{(k)}_{i+1} &= \textrm{prox}_{\bar{B}_{\mu}/\beta}(\tilde{\bu}^{(k)}_{i+1} - \bv^{(k)}_{i}), \label{cp-prox} \\
	\bv^{(k)}_{i+1} &= \bv^{(k)}_{i} - \tilde{\bu}^{(k)}_{i+1} + \bu^{(k)}_{i+1}, \label{cp-dual}
\end{align}
where $\bar{B}_{\mu}(\cdot)$ is the log-barrier function $B(\bu,\bv,\mu)$ restricted to the $\bu$ part, that is, 
\[
\bar{B}_{\mu}(\bu) =\beta((\bx^0)^{T} \bs^0+1)\theta+\mathds{1}(\xi=-(\bx^0)^{T} \bs^0-1)+\mu F(\bx)-\mu\log\tau
\]
The proximal operator is defined by $\prox_f(\bx)=\argmin_{\bar{\bu}}\{f(\bar{\bu})+\frac{1}{2}\|\bar{\bu}-\bu\|^2\}$.

It is clear that, ABIP needs to the inverse of $I+Q$, which can be computed before the algorithm proceed. The remaining computational bottleneck is to compute the proximal operator. In Section \ref{subsec:prox-cone}, we show that the proximal operator associated with several widely-used convex cones can be solved efficiently.

\paragraph{Termination criteria.} For the inner loop, we run ADMM iteration until the following termination criterion is satisfied: 
\begin{equation}
    ||\bQ\bu_{i}^{(k)}-\bv_{i}^{(k)}||^2\leq \mu^{(k)}.\label{inner-loop-termination}
\end{equation}
For the outer loop, we theoretically terminate the algorithm when $\mu^{(k)} < \epsilon$. In the implementation, however, a different termination criterion is employed to align with other solvers, which will be introduced later.

\subsection{Convergence analysis}
In this subsection, we develop the iteration complexity of ABIP, for which we generalize the analysis of \cite{ABIP-2021} to linear conic optimization. To begin with, we outline the key ingredients of the analysis.
First, we prove that the iterates generated by ADMM in each inner loop are uniformly bounded above. 
Then, based on this property, we show that ABIP exhibits a linear convergence to the optimal solution associated with each subproblem. 
Finally, we derive the total iteration complexity of ABIP under the path-following framework.

Let $(\bu_k^*,\bv_k^*)$ be the optimal solution to subproblem~\eqref{conic-subproblem}. In the following proposition, we show that the iterates generated by ADMM will converge to $(\bu_k^*,\bv_k^*)$. The proof can be extended from \cite{ABIP-2021}, and thus is omitted in the paper.
\begin{proposition}
\label{prop:inner-converge}
For a fixed $k\in\mathbb{N}_+$, the sequence $\{||\bu_i^{(k)}-\bu_k^*||^2+||\bv_i^{(k)}-\bv_k^*||^2\}_{i\ge 0}$ is monotonically decreasing and converges to 0.\label{prop:monotone-decrease}
\end{proposition}

In the next proposition, we further prove that the sequence $\{||\bu_i^{(k)}-\bu_k^*||^2+||\bv_i^{(k)}-\bv_k^*||^2\}_{i\ge 0}$ is uniformly bounded. With Proposition \ref{prop:inner-converge}, it suffices to show that the initial point $\{||\bu_0^{(k)}-\bu_k^*||^2+||\bv_0^{(k)}-\bv_k^*||^2\}$ is bounded.
\begin{proposition}\label{prop:bounded-seq}
    The sequence $\{||\bu_i^{(k)}-\bu_k^*||^2+||\bv_i^{(k)}-\bv_k^*||^2\}$ is uniformly bounded.
\end{proposition}

\proof{Proof.}
We first recall an important fact that the set of central path points $\{(\bu_k^*, \bv_k^*)\}_{k\geq 1}$ is uniformly bounded. That is, there exists a constant $C_1$ such that 
\begin{equation*}
    ||\bu_k^*||^2+||\bv_k^*||^2\leq C_1.
\end{equation*}

Let $(\bu^*,\bv^*)$ be the limit point of $\{(\bu_k^*, \bv_k^*)\}_{k\geq 1}$, then $(\bu^*,\bv^*)$ also satisfies $||\bu^*||^2+||\bv^*||^2\leq C_1$.
Let $N_k$ denote the number of inner iterations in the $k$-th outer loop. With Proposition \ref{prop:inner-converge} and the inner loop termination criterion \eqref{inner-loop-termination}, we have $N_k < \infty$ for any fixed $k$. 

Now we claim that when $k\to\infty$, $||\bu_{N_k}^{(k)}-\bu^*||^2+||\bv_{N_k}^{(k)}-\bv^*||\to 0$. Otherwise, there exists $\delta>0$ and a subsequence $\{k_t\}_{t\in\mathbb{N}^+}$ with $k_t \uparrow \infty$ as $t \rightarrow \infty$ such that $||\bu_{N_{k_t}}^{(k_t)}-\bu^*||^2+||\bv_{N_{k_t}}^{(k_t)}-\bv^*||>\delta$ for all $t\in\mathbb{N}^+$, and there exists an integer $T$ such that $||\bu_{N_{k_t}}^{(k_t)}-\bu_{k_t}^*||^2+||\bv_{N_{k_t}}^{(k_t)}-\bv_{k_t}^*||>\delta/2$ for all $t>T$, which contradicts with~\eqref{inner-loop-termination}. Therefore, $||\bu_{N_k}^{(k)}-\bu^*||^2+||\bv_{N_k}^{(k)}-\bv^*||\to 0$ when $k\to\infty$, which also implies that there exists $C_2 > 0 $ such that
\begin{equation*}
||\bu_{N_k}^{(k)}||^2+||\bv_{N_k}^{(k)}||^2 \leq C_2.
\end{equation*}

Since $\bu_{N_k}^{(k)}=\bu_{0}^{(k+1)}$ and $\bv_{N_k}^{(k)}=\bv_{0}^{(k+1)}$, we have proven the boundedness of the sequence $\{\bv_0^{(k)}\}$. Moreover, by the results of Proposition~\ref{prop:monotone-decrease}, we have that $\{(\bu_i^{(k)},\bv_i^{(k)})\}$ is uniformly bounded.
 \hfill $\square$
\endproof

With the above results, we give an upper bound on the number of ADMM iterations for each inner loop in the next proposition, and present the total iteration complexity of ABIP in Theorem \ref{thm:total-bound}.
\begin{proposition}
\label{prop:complexity-admm}
The number of ADMM iterations in each inner loop associated with $\mu^{(k)}$, denoted by $N_k$, satisfies
	\begin{equation*}
		N_k\leq\log\left(\dfrac{4C_1(1+||\bQ||^2)}{\mu^{(k)}}\right)\left[\log\left(1+\min\left\{\dfrac{1}{C_3},\dfrac{\mu^{(k)}}{2C_DC_3\beta}\right\}\right)\right]^{-1},
	\end{equation*}
	where
	\begin{equation*}
		C_3=\left[1+\dfrac{12\lambda_{\max}(\bA^{T}\bA)}{\lambda^2_{\min} (\bA\bA^{T})}\max\{1,||\bc||^2,||\br_d||^2\}\right]\cdot \left[1+\dfrac{6}{||\br_p||^2}\max\{||\bb||^2,||\bA||^2\}\right].
	\end{equation*}
 and $1/C_D > 0$ is the strong convexity parameter of the log-barrier function.
\end{proposition}
The proof can be found in Appendix \ref{proof-prop-5}. Now we are ready to develop the total iteration complexity of ABIP for solving linear conic optimization problems.
\begin{theorem}
\label{thm:total-bound}
Suppose that ABIP is terminated when $\mu^{(k)}<\epsilon$, where $\epsilon > 0$ is a pre-specified tolerance. It requires a total number of $T_1=\mathcal{O}\left(\log(1/\epsilon)\right)$ outer loops and $T_2=\mathcal{O} \left( \kappa_A^2||\bQ||^2 / \epsilon \cdot \log \left( 1 / \epsilon \right) \right)$ ADMM iterations, where $\kappa_{\bA}:=\lambda_{\max}(\bA^T\bA)/\lambda_{\min}(\bA \bA^T)$.
\end{theorem}
\proof{Proof.}
Note that ABIP is a double-loop algorithm. The outer loop is terminated when $\mu^{(k)}<\epsilon$, where $\epsilon>0$ is a pre-specified tolerance level. It is easy to see that the number of outer iterations is
\begin{equation*}
T_{\texttt{IPM}}=\Big\lceil\dfrac{\log(\mu^{(0)}/\epsilon)}{\log(1/\phi)}\Big\rceil.
\end{equation*}
For the total number of ADMM iterations, we have the following bound:
\begin{align*}
    T_{\texttt{ADMM}}&=\sum_{k=1}^{T_{\texttt{IPM}}}N_k\leq\sum_{k=1}^{T_{\texttt{IPM}}}\log\left(\dfrac{4C_1(1+||\bQ||^2)}{\mu^{(k)}}\right)\left[\log\left(1+\min\left\{\dfrac{1}{C_3},\dfrac{\mu^{(k)}}{2C_DC_3\beta}\right\}\right)\right]^{-1}\\
    &=\sum_{k=1}^{T_{\texttt{IPM}}}\log\left(\dfrac{4C_1(1+||\bQ||^2)}{\mu^{(0)}\phi^{k}}\right)\left[\log\left(1+\min\left\{\dfrac{1}{C_3},\dfrac{\mu^{(k)}}{2C_DC_3\beta}\right\}\right)\right]^{-1}\\
    &=\mathcal{O}\left(\dfrac{\kappa_A^2||\bQ||^2}{\epsilon}\log\left(\dfrac{1}{\epsilon}\right)\right).
\end{align*}

\hfill $\square$
\endproof

With Theorem \ref{thm:total-bound}, we now make a comparison of iteration complexity between ABIP and PDLP in Table \ref{tab:complexity}. It shows that ABIP is more sensitive to the condition number of the coefficient matrix compared with PDLP, which also motivates us to improve the practical performance of ABIP in the next section. It is important to note, however, that our result is established for general linear conic optimization, whereas PDLP's analysis is specifically tailored for linear optimization.

To conclude this subsection, we further give the total arithmetic operations of ABIP, that is, $\mathcal{O} \left(n^3 + \left( n^2 \kappa_{\bA}^2 \|\bQ\|^2 / \epsilon \right) \log{(1/\epsilon)} \right) $. The first term $\mathcal{O}(n^3)$ comes from the matrix decomposition, and the second term is from the $\mathcal{O}(n^2)$ arithmetic operations required by each ADMM iteration and the total number of ADMM iterations, that is, $ \mathcal{O}\left(( \kappa_{\bA}^2 \|\bQ\|^2/\epsilon) \log{({1}/{\epsilon})} \right)$.
\begin{table}[htbp]
    \centering
    \begin{tabular}{c|c|c}
    \hline
    Algorithm & Problem &  Complexity  \\
    \hline
    \rule{0pt}{3ex} ABIP  & Conic Optimization &  $\mathcal{O} \left( \kappa_A^2||\bQ||^2 / \epsilon \cdot \log \left( 1 / \epsilon \right) \right)$     \\
    \rule[-1.5ex]{0pt}{0pt} PDLP \citep{PDHG-2021} & Linear Optimization &  $\mathcal{O}\left( \sigma_{\max}(\bA) / \sigma_{\min}^{+}(\bA)  \cdot \log\left(1/\epsilon\right) \right)$ 
    \\
     \hline 
    \end{tabular}
    \caption{Iteration Complexity of ABIP and PDLP. Note that $\kappa_{\bA}$, $\sigma_{\min}^{+}(\bA)$ and $\sigma_{\max}(\bA)$ represent the condition number of $\bA$, the minimum nonzero singular value of $\bA$, and the maximum singular value of $\bA$, respectively. The matrix $\bQ$ is the coefficient matrix of the HSD formulation in \eqref{hsd2}.}
    \label{tab:complexity}
\end{table}

\subsection{Solving the linear system}
In this subsection, we discuss how to efficiently solve the linear system of $(\bI + \bQ) \tilde{\bu} = \bw$ for some given $\bw$. In particular, we want to solve
\begin{equation} \label{eq:linsys}
\begin{bmatrix}
        \bM & \bh \\
        -\bh^T& 1
    \end{bmatrix} 
    \tilde{\bu}=\bw,
\end{equation} 
where
\begin{equation*}
    \bM = \begin{bmatrix}
        \bI_{m} & \bA\\
        -\bA^T & \bI_{n} 
    \end{bmatrix} \quad \text{and} \quad  \bh = \begin{bmatrix}
        -\bb \\
        \bc
    \end{bmatrix}.
\end{equation*}
To further simplify \eqref{eq:linsys}, one can follow \cite{ABIP-2021} and only consider a smaller linear system restricted to the $(\tilde{\by}, \tilde{\bx})$ block variables of $\tilde{\bu}$:
\begin{equation}\label{eq:linsys-reduced}
\begin{bmatrix}
\tilde{\by}\\
\tilde{\bx}
\end{bmatrix}
=\bM^{-1}
\begin{bmatrix}
\by+\br\\
\bx +\bs
\end{bmatrix}.
\end{equation}

To solve this, the first approach, or the so-called direct method, is to perform a sparse permuted $\bL \bD \bL^{\top}$ factorization of $\bM$ before the first iteration, and store the factors $\bL$ and $\bD$. Note that the factorization only needs to be computed one time, and the subsequent ADMM iterations can be carried out by solving much easier triangular and diagonal linear systems with the cached factors $\bL$ and $\bD$. The second approach, or the indirect method, is to approximately solve the system with conjugate gradient method \citep{wright1999numerical}. When the time or memory cost of factorizing $\bM$ is expensive, the indirect method is more favorable due to its simple and efficient update rule. Note that similar strategies have been employed in many other popular open-source solvers, such as OSQP~\citep{stellato2020osqp} and SCS~\citep{SCS-2016}.
 
Although the general-purpose approaches designed for solving linear systems are convenient, they may still have limitations since the inner structure of the problem is hardly exploited. 
This further prohibits the use of ABIP for large-scale problems in machine learning and data science. 
Therefore, to accelerate the performance of ABIP for such applications, it is desirable to use a more specialized linear system solver that can be adaptive to the problem structure. 
For some important problems in machine learning, such as LASSO and SVM, we present novel implementation strategies to solve linear systems, which significantly improve the scalability of ABIP for large-scale problems. 
The more details are presented in Appendix~\ref{sec:linsys}.

\subsection{Solving the proximal problem}\label{subsec:prox-cone}
In this subsection, we discuss how to solve the proximal problem in \eqref{cp-prox}. \cite{ABIP-2021} have shown that, the proximal problem associated with the nonnegative orthant exhibits a closed-form solution. As a complement, we show that the proximal problems associated with other convex cones, including second-order cone and semidefinite cone, can also be computed efficiently.

\subsubsection{Second-order cone}
We consider the following proximal problem:
\begin{equation}
	\bx^{*} = \argmin_{\bx \in \text{int}\mathcal{K}} \  \lambda F(\bx) + \frac{1}{2}\|\bx - \bzeta\|_{2}^{2}, \label{subp}
\end{equation}
where $\lambda > 0$ and $\bzeta$ is a column vector. In this part, we show that the proximal problems associated with standard second-order cone (SOC) and rotated second-order cone (RSOC) have closed-form solutions. We present the main results here and leave the technical details to Appendix \ref{append-prox-soc}.

\paragraph{Proximal problem with standard SOC.} For $\bx = (t, \bar{\bx}) \in \text{SOC}$, the log-barrier function is defined by $F(\bx) = - \log \left( t^{2} - \bar{\bx}^{T}\bar{\bx} \right)$. For convenience, we let $\bzeta = ( \zeta_{t} , \bzeta_{\bar{\bx}} )^{T}$ and discuss the following two situations. When $ \zeta_{t} = 0 $, we have
\begin{equation*}
    t^{*} = \sqrt{2\lambda + \frac{1}{4}\bzeta_{\bar{\bx}}^{T}\bzeta_{\bar{\bx}}}, \quad
    \bar{\bx}^{*} = \frac{1}{2}\bzeta_{\bar{\bx}}.
\end{equation*}
When $ \zeta_{t} \neq 0 $, we define $\rho_{1}$ and $\rho_{2}$, respectively, by
\begin{equation*}
	\rho_{1} = \frac{\gamma - \sqrt{\gamma^{2}-16}}{2} < 2, \quad  \rho_{2} = \frac{\gamma + \sqrt{\gamma^{2}-16}}{2} > 2,
\end{equation*}
where
\[
\gamma = \frac{\frac{\zeta_{t}^{2} - \bzeta_{\bar{\bx}}^{T}\bzeta_{\bar{\bx}}}{\lambda} + \sqrt{\left( \frac{\zeta_{t}^{2} - \bzeta_{\bar{\bx}}^{T}\bzeta_{\bar{\bx}}}{\lambda}\right)^{2} + 4\left( \frac{4\zeta_{t}^{2}+ 4\bzeta_{\bar{\bx}}^{T}\bzeta_{\bar{\bx}}}{\lambda} + 16 \right) }}{2}.
\]
Then, we can obtain the solution
\begin{equation*}
	t^{*} = \frac{\rho}{\rho - 2}\zeta_{t}, \quad
	\bar{\bx}^{*} = \frac{\rho}{\rho+2}\bzeta_{\bar{\bx}},
\end{equation*}
where $\rho=\rho_{1}$ if $\zeta_{t} < 0$ and $\rho=\rho_{2}$ otherwise.

\paragraph{Proximal problem with standard RSOC.} Let $\bx = (\eta,\nu,\bar{\bx}) \in \text{RSOC}$ and define the log-barrier function by $F(\bx) = -\log \left(\eta\nu - \frac{1}{2}\bar{\bx}^{T}\bar{\bx} \right)$. For convenience, we let $\bzeta = ( \zeta_{\eta} , \zeta_{\nu} , \bzeta_{\bar{\bx}} )^{T}$ and discuss the following two situations. If $\zeta_{\eta} + \zeta_{\nu} = 0$, we have
\begin{equation*}
        \eta^{*} = \frac{\zeta_{\eta} + \sqrt{\zeta_{\eta}^{2} + 4\left(\lambda + \frac{1}{8}\bzeta_{\bar{\bx}}^{T}\bzeta_{\bar{\bx}} \right) }}{2}, \quad
		\nu^{*} = \frac{-\zeta_{\eta} + \sqrt{\zeta_{\eta}^{2} + 4\left(\lambda + \frac{1}{8}\bzeta_{\bar{\bx}}^{T}\bzeta_{\bar{\bx}} \right) }}{2}, \quad
        \bar{\bx}^{*} = \frac{1}{2}\bzeta_{\bar{\bx}}.
\end{equation*}
When $\zeta_{\eta} + \zeta_{\nu} \neq 0$, we define $\rho_{1}$ and $\rho_{2}$ as
\begin{equation*}
	\rho_{1} = \frac{\gamma - \sqrt{\gamma^{2}-4}}{2} < 1, \quad \rho_{2} = \frac{\gamma + \sqrt{\gamma^{2}-4}}{2} > 1,
\end{equation*}
where
\begin{equation*}
	\gamma = \frac{\frac{2\zeta_{\eta}\zeta_{\nu} - \bzeta_{\bar{\bx}}^{T}\bzeta_{\bar{\bx}}}{2\lambda} + \sqrt{\left( \frac{2\zeta_{\eta}\zeta_{\nu} - \bzeta_{\bar{\bx}}^{T}\bzeta_{\bar{\bx}}}{2\lambda}\right)^{2} + 4\left( \frac{\zeta_{\eta}^{2} + \zeta_{\nu}^{2} + \bzeta_{\bar{\bx}}^{T}\bzeta_{\bar{\bx}}}{\lambda} + 4 \right) }}{2}.
\end{equation*}
Then, we have
\begin{equation*}
	\eta^{*} = \frac{\rho^{2}\zeta_{\eta} + \rho\zeta_{\nu}}{(\rho+1)(\rho-1)}, \quad
	\nu^{*} = \frac{\rho\zeta_{\eta} + \rho^{2}\zeta_{\nu}}{(\rho+1)(\rho-1)} \quad
	\text{and}\quad \bar{\bx}^{*} = \frac{\rho\bzeta_{\bar{\bx}}}{\rho+1},
\end{equation*}
where $\rho = \rho_{1}$ if $\zeta_{\eta} + \zeta_{\nu} < 0$ and $\rho = \rho_{2}$ otherwise.

\subsubsection{Semidefinite cone}
In this part, we consider the proximal problem associated with the semidefinite cone (SDC), which can be represented by $\mathcal{K}=\{\bX: \bX\in\mathbb{R}^{n\times n}, \bX\succeq 0, \bX=\bX^T \}$. 
The log-barrier function is defined by $F(\bX)=-\log\det(\bX)$. We consider the following proximal subproblem:
\begin{equation}
    \bX^*=\argmin_{\bX \succ 0, \bX=\bX^T} -\lambda \log\det(\bX)+\dfrac{1}{2}||\bX-\bA||^2,\label{proximal-sdc}
\end{equation}
where $\lambda > 0$, $\bA\in\mathbb{R}^{n\times n}$ is a symmetric matrix, and $||\cdot||$ represents matrix Frobenius norm, i.e., $||\bA||=\sqrt{\text{tr}(\bA^{T} \bA)}$. By the optimality condition, problem \eqref{proximal-sdc} is equivalent to finding the root of the following equation:
\begin{equation}
\label{eq:sdp-equation}
    -\lambda (\bX)^{-1}+\bX-\bA=\bzero,\ \bX\succ 0.
\end{equation}

Since $\bA$ is real symmetric, we have its spectral decomposition $\bA=\bQ^T \bD\bQ$, where $\bD=\text{diag}(d_i)_{i=1}^n$ is a diagonal matrix with eigenvalues $d_i$, and $\bQ$ is an orthogonal matrix. Then, we define a diagonal matrix
\begin{equation*}
   \bE=\text{diag}(e_i)_{i=1}^n, \text{ where } e_i=\dfrac{d_i+\sqrt{d_i^2+4\lambda}}{2},\ \forall i  = 1, 2, \ldots, n.
\end{equation*}
Moreover, it is easy to verify that
\begin{equation*}
    \lambda \bE^{-1}+\bE-\bD=\bzero,\ \bE\succ 0,
\end{equation*}
which further implies
\begin{equation*}
    \bQ^{T}(\lambda \bE^{-1}+\bE-\bD)\bQ = \lambda \bQ^T \bE^{-1}\bQ+\bQ^{T} \bE\bQ-\bA=\bzero.
\end{equation*}
In other words, we have shown that $\bQ^{T}\bE\bQ$ is a solution to equation \eqref{eq:sdp-equation}. Due to the strong convexity, the optimal solution $\bX^*$ to the proximal problem~\eqref{proximal-sdc} must be unique, which is exactly given by $\bX^*=\bQ^T \bE\bQ$.

\section{Newly implemented strategies of ABIP for linear optimization}
\label{Sec:new-strategy}

In this section, we introduce several new strategies that further improve ABIP's performance for LP. Specifically, the LP problem has the following form: 
\begin{equation}\label{eq:lp}
    \begin{split}
        \min & \quad  {\bc}^{T} {\bx} \\
        \text{s.t.} & \quad {\bA} {\bx} = {\bb}, \\
                    &\quad {\bx} \geq \bzero. \\
    \end{split}
\end{equation}
We remark that, although these strategies are applicable to general conic optimization, our preliminary experiments show that their benefits are considerably more pronounced for LP. Consequently, our subsequent discussion will focus on the enhancements within the realm of LP. A comprehensive empirical study is deferred to Section \ref{Sec:numerical}.

\subsection{An adaptive strategy for the barrier parameter $\mu$}
Under the framework of ABIP, the barrier parameter $\mu^{(k)}$ plays a critical role in
balancing the subproblem's convexity and global optimality of the obtained ADMM iterates. 
Moreover, a proper path driving $\mu^{(k)}$ down to $0$ benefits the convergence of ABIP by reducing the total number of both IPM iterations and ADMM iterations. In particular, we consider the following two popular adaptive methods for the barrier parameter update.

\paragraph{The aggressive strategy \citep{wachter2006implementation}.} Given  $0<\zeta<1$ and $\eta>1$, the aggressive strategy updates by 
\[ \mu^{(k + 1)} = \min \{ \zeta \mu^{(k)}, (\mu^{(k)})^{\eta} \} . \]
This is often used at the beginning of the algorithm when the subproblem is well-behaved and allows an aggressive decrease in the barrier parameter. However, as ABIP progresses, this aggressive strategy will result in more but useless ADMM iterations to solve the subproblem since $\mu$ will be quite small. Hence, under this case, we will switch to the following strategy to gradually decrease $\mu$. 

\paragraph{The LOQO strategy \citep{vanderbei1999loqo}.} Given parameter $0<\alpha<1$, the LOQO strategy updates by 
\[ \mu^{(k + 1)} = \mu^{(k)} \cdot \max \left\{ 0.1 \cdot \min \left\{ 0.05
   \cdot \frac{1 - \phi}{\phi}, 2 \right\}^3, \alpha \right\}, \]
where $\phi := (n + 2) \cdummy \min \{ \min_{0 \leq i \leq n} \{ x_i s_i \}, \tau
\kappa, \theta \xi \} / (x^T s + \kappa \tau + \theta \xi)$. This method is
adapted from the well-known update rule in the quadratic optimization solver
LOQO, where $\phi$ serves as a measure of centrality and $\mu^{(k + 1)}$ is
lower-bounded by $\alpha \mu^{(k)}$ to avoid over-aggressive updates.

In our implementation, the two methods are integrated into a hybrid strategy that
automatically switches between them. To be more specific, we switch from the aggressive strategy to the LOQO strategy when the barrier parameter $\mu^{(k)} \leq 10^3 \cdot \epsilon $, where $\epsilon$ is the pre-determined tolerance level. %

\subsection{Restart scheme}
The restart scheme is a useful technique to improve the
practical and theoretical performance of an optimization algorithm \citep{Pokutta-2020}.
Motivated by the acceleration performance of the restart scheme for solving large-scale linear optimization \citep{PDLP-2021, PDHG-2021}, we implement a restart scheme in the inner loop of ABIP.

In particular, let $M$ be the accumulated number of ADMM iterations performed by ABIP so far, and $M_{k}$ be the number of ADMM iterations for solving the inner problem \eqref{conic-subproblem} at the $k$-th outer loop.
We observe that the restart scheme tends to be less effective in the initial stage of ABIP. Therefore, we employ this strategy only when the accumulated number of iterations $M$ is sufficiently large, such as surpassing a pre-determined threshold $M_T$. 
Once this scheme is triggered, the ADMM algorithm in the inner loop will be restarted with a fixed period $F$.
Namely, we restart the algorithm if
\begin{equation*}
    M \geq M_T \quad  \text{and} \quad (M_{k}+1) \bmod F = 0.
\end{equation*}
Then, we use the average of historical iterates in the current restart cycle as the initial point of the next restart cycle.
Details of the restart scheme are presented in Algorithm \ref{alg:abip-with-restart}.

\begin{algorithm}[ht]
\SetAlgoNoLine
\caption{ABIP with fixed frequency restart}
\label{alg:abip-with-restart}
Initialize $\mathbf{u}^{(0)}_{0}$ and $\mathbf{v}^{(0)}_{0}$ \;
Set ADMM iteration counter $M \leftarrow 0$ \;
Set restart threshold $M_{T}$ and restart period $F$ \;
\For{$k=0, 1, 2,\cdots$}{
    $\hat{\mathbf{u}}^{(k)}_{\text{avg}} \leftarrow \mathbf{u}^{(k)}_{0} $, $\hat{\mathbf{v}}^{(k)}_{\text{avg}} \leftarrow \mathbf{v}^{(k)}_{0} $ \;
    $M_k \leftarrow 0$ \;
    \While{True}{
        \If{the \textbf{inner} termination condition is satisfied}{
            break \;}
        $\mathbf{u}_{M_{k}+1}^{(k)}$, $\mathbf{v}_{M_{k}+1}^{(k)} \leftarrow $ ADMM Step$(\mathbf{u}_{M_{k}}^{(k)}$, $\mathbf{v}_{M_{k}}^{(k)})$ \;
        $\hat{\mathbf{u}}^{(k)}_{\text{avg}} \leftarrow \hat{\mathbf{u}}^{(k)}_{\text{avg}} + \mathbf{u}_{M_{k}+1}^{(k)}$ \;
        $\hat{\mathbf{v}}^{(k)}_{\text{avg}} \leftarrow \hat{\mathbf{v}}^{(k)}_{\text{avg}} + \mathbf{v}_{M_{k}+1}^{(k)}$ \;
        $M_k \leftarrow M_k + 1$ \;
        $M \leftarrow M+1$ \;
        \If{$M\geq M_T$ and $ M_{k} \bmod F = 0$}{
        $\mathbf{u}_{M_{k}+1}^{(k)} \leftarrow \hat{\mathbf{u}}_{\text{avg}}^{(k)} / F $, $\mathbf{v}_{M_{k}+1}^{(k)} \leftarrow \hat{\mathbf{v}}_{\text{avg}}^{(k)} / F $ \;
        $ \hat{\mathbf{u}}_{\text{avg}}^{(k)} \leftarrow \mathbf{u}_{M_{k}+1}^{(k)} $, $\hat{\mathbf{v}}_{\text{avg}}^{(k)} \leftarrow \mathbf{v}_{M_{k}+1}^{(k)} $ \;
        }
        \If{the \textbf{final} termination criterion is satisfied}{
        \Return{}\;}
    }
$\mathbf{u}^{(k+1)}_{0} \leftarrow \mathbf{u}^{(k)}_{M_{k}} $, $\mathbf{v}^{(k+1)}_{0} \leftarrow \mathbf{v}^{(k)}_{M_{k}} $  \;    
}
\end{algorithm}

\begin{figure}[h]
    \centering
    \subfigure[The spiral trajectory of the original ABIP]{
    \begin{minipage}[t]{0.40\linewidth}
    \centering
    \includegraphics[scale=0.5]{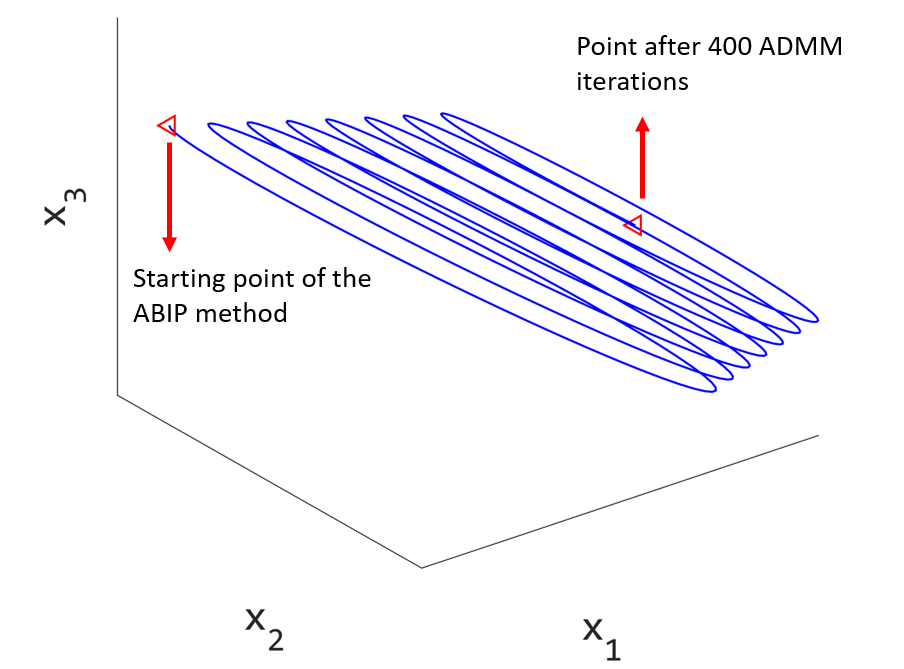}
    \end{minipage}
    }
    \hfill
    \subfigure[The trajectory of ABIP after restart]{
    \begin{minipage}[t]{0.40\linewidth}
    \centering
    \includegraphics[scale=0.48]{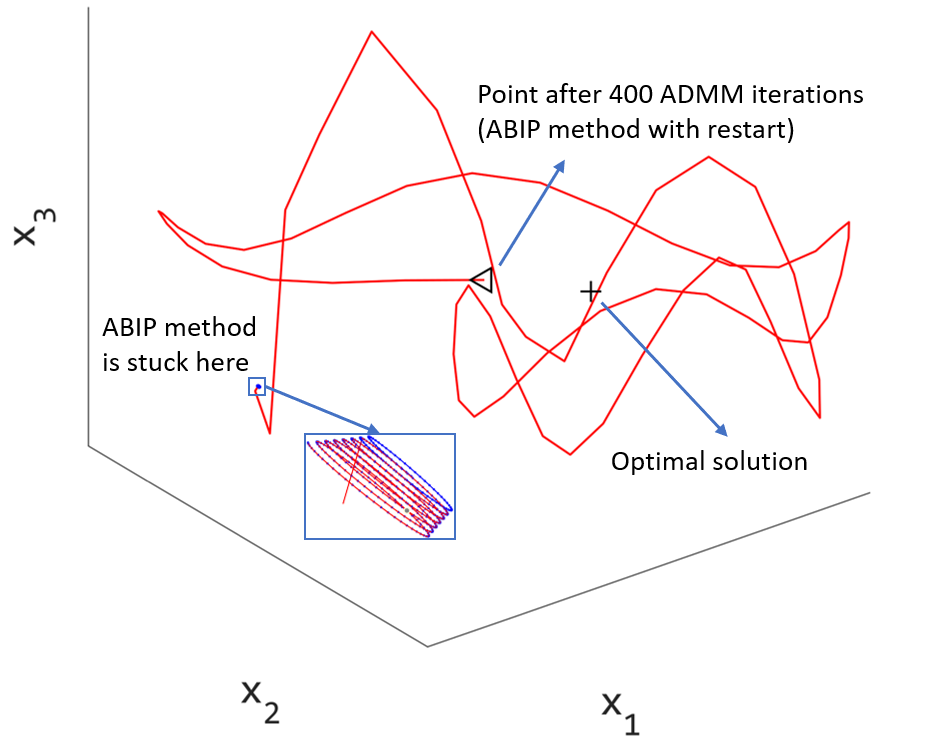}
    \end{minipage}
    }
    \caption{Comparison of the trajectories of ABIP after restart and the origin ABIP}
    \label{fig:restart}
\end{figure}

To illustrate the acceleration performance of the restart strategy, we visualize the trajectory of ABIP for solving the \texttt{SC50B} instance from Netlib as an example in Figure \ref{fig:restart}.
For ease of exposition, we limit our visualization to the first three dimensions of the candidate solutions.
Specifically, Figure \ref{fig:restart}(a) shows that the original ABIP has a spiral trajectory, leading to a long time to converge to an optimal solution. 
However, ABIP progresses more aggressively after restart and is closer to the optimal solution than the original ABIP due to averaging the previous iterates, as shown in Figure \ref{fig:restart}(b). 
In fact, the restart strategy contributes to reducing almost $70\%$ total ADMM iterations on the \texttt{SC50B} instance; see Table \ref{tab:restart-example}. Preliminary experiments also show that the restart strategy is more efficient when one hopes to obtain a high-accuracy solution, for example, $\epsilon = 10^{-6}$.

\begin{table}[h]
    \centering
    \begin{tabular}{c|c|c}
    \hline
     &  ABIP  & ABIP with restart  \\
     \hline
    \# of ADMM iterations &  309175  & \textbf{93313} \\ 
    \hline
    \end{tabular}
    \caption{\# of ADMM iterations for \texttt{SC50B} of ABIP and ABIP with restart strategy, $\epsilon = 10^{-6}$, direct method}
    \label{tab:restart-example}
\end{table}

\subsection{New inner loop termination criteria}

The original ABIP adopts the criterion \eqref{inner-loop-termination} to terminate the inner loop. 
In its enhanced implementation,
motivated by the restart strategy, we also compute the average of previous iterates $ \bar{\mathbf{u}}_{i}^{(k)} = \sum_{i = 1}^{k} \mathbf{u}_{i}^{(k)} / k $ and $\bar{\mathbf{v}}_{i}^{(k)} =  \sum_{i = 1}^{k} \mathbf{v}_{i}^{(k)} / k $, and adopt $\left\|\bQ \bar{\mathbf{u}}_{i}^{(k)}-\bar{\mathbf{v}}_{i}^{(k)}\right\|^{2} \leq \mu^{(k)}$ as another stopping criterion when solving the $k$-th subproblem. When one of the two criteria is satisfied, we terminate the inner loop and let the returned solution be the initial point of ADMM to solve the subproblem in the form of \eqref{conic-subproblem} corresponding to $\mu^{(k+1)}$.

As mentioned before, we use a different criterion in the implementation to check the global convergence. Specifically, we stop the algorithm when the following condition is satisfied:
\begin{equation}
\label{eq:three-all-small}
    \max \left(\epsilon_{\texttt{pres}}^{\prime}, \epsilon_{\texttt{dres}}^{\prime}, \epsilon_{\texttt{dgap}}^{\prime}  \right) \leq \epsilon,
\end{equation}
where $ \epsilon_{\texttt{pres}}^{\prime}, \epsilon_{\texttt{dres}}^{\prime}$, and $\epsilon_{\texttt{dgap}}^{\prime}$ represent the primal residual, dual residual, and duality gap, respectively. For linear optimization, we compute $ \epsilon_{\texttt{pres}}^{\prime}, \epsilon_{\texttt{dres}}^{\prime}, \epsilon_{\texttt{dgap}}^{\prime} $ by 
\begin{equation}\label{outer:criteria}
    \begin{aligned}
        \epsilon_{\texttt{pres}}^{\prime} &= \left\|\bA\frac{\bx}{\tau} - \bb\right\| / \left(1 + \|\bb\|\right),   \\
        \epsilon_{\texttt{dres}}^{\prime} &= \left\|\bA^T\frac{\by}{\tau} + \frac{\bs}{\tau} - \bc\right\| /  \left(1 + \|\bc\| \right), \\
        \epsilon_{\texttt{dgap}}^{\prime} &= \left|\bc^T\frac{\bx}{\tau} - \bb^T\frac{\by}{\tau}\right| /  \left(1 + \left|\bc^T\frac{\bx}{\tau}\right| + \left|\bb^T\frac{\by}{\tau}\right|\right).
    \end{aligned}
\end{equation}
In the inner loop, we use \eqref{outer:criteria} to check the global convergence. This is an early-stop strategy to jump out of the inner loop, which is particularly useful when some ADMM iterates do not meet the inner loop convergence criterion but already satisfy the global convergence criterion. We remark that when to activate this strategy should be chosen carefully. If we check the global convergence criterion too early, the algorithm will waste a considerable amount of time to verify those ADMM iterates that are far from the global optimality. In the implementation, we employ this strategy when the barrier parameter $\mu^{(k)}$ satisfies $\mu^{(k)} < \epsilon$. Although ABIP will terminate in theory when such condition is satisfied, the obtained iterate does not satisfy \eqref{outer:criteria} in practice, and thus the algorithm still proceeds.

\subsection{Half update}

As mentioned before, \cite{ABIP-2021} propose a customized initialization of $(\bp_0^{(k)}, \bq_0^{(k)})$ for the inner problem \eqref{conic-subproblem} to eliminate the variables $\bp,\bq$ and $\Tilde{\bv}$ in the update rule, and then sequentially update $\Tilde{\bu}^{(k)}_{i+1}, \bu^{(k)}_{i+1}, \bv^{(k)}_{i+1}$ in each ADMM iteration. 
Here, $\bv^{(k)}_{i+1}$ can be seen as the role of dual variable. 
This further motivates us to introduce the half update strategy in ADMM-type algorithm by updating $\bv^{(k)}$ twice.
In particular, we first update the primal variable $\tilde{\mathbf{u}}_{i+1}^{(k)}$, followed by the update of dual variable $\mathbf{v}_{i+\frac{1}{2}}^{(k)}$ with stepsize $\alpha_1\in(0,2)$. Then, we update the primal variable $\mathbf{u}_{i+1}^{(k)}$, followed by the update of dual variable $\mathbf{v}_{i+1}^{(k)}$ with stepsize $\alpha_2\in(0,2)$. 
Therefore, the update rule for $\mathbf{u}^{(k)}_{i+1}$ is slightly changed as
\begin{equation}
\label{eq:update-mu}
\mathbf{u}_{i+1}^{(k)}=\underset{\mathbf{u}}{\operatorname{argmin}}\left[\bar{B}\left(\mathbf{u}, \mu^{(k)}\right)+\frac{\beta}{2}\left\|\mathbf{u}-\tilde{\mathbf{u}}_{i+1}^{(k)}+ \mathbf{v}_{i+\frac{1}{2}}^{(k)}\right\|^{2}\right],
\end{equation}
where $\bar{B}(\mathbf{u}, \mu^{(k)}) = \beta(n+1) \theta-\mu^{(k)} \log (\mathbf{x})-\mu^{(k)} \log \tau $.
Meanwhile, problem \eqref{eq:update-mu} admits closed-form solutions given by
\begin{align}
    \mathbf{y}_{i+1}^{(k)} &= \arg\min_{\mathbf{y}} \frac{\beta}{2} \| \mathbf{y} - \tilde{\mathbf{y}}^{(k)}_{i+1} + \mathbf{r}^{(k)}_{i+\frac{1}{2}} \|^{2} = \tilde{\mathbf{y}}^{(k)}_{i+1} - \mathbf{r}^{(k)}_{i+\frac{1}{2}} \label{update:y} \\
    \mathbf{x}_{i+1}^{(k)} &= \arg\min_{\mathbf{x}}  \left[ - \frac{\mu^{(k)}}{\beta} \log{(\mathbf{x})} + \frac{1}{2} \left\| \mathbf{x} - \tilde{\mathbf{x}}^{(k)}_{i+1} + \mathbf{s}^{(k)}_{i+\frac{1}{2}}  \right\|^{2}  \right] \nonumber \\
    &= \frac{1}{2}\left[\left(\tilde{\mathbf{x}}_{i+1}^{(k)}-\mathbf{s}_{i+\frac{1}{2}}^{(k)}\right)+\sqrt{\left(\tilde{\mathbf{x}}_{i+1}^{(k)}-\mathbf{s}_{i+\frac{1}{2}}^{(k)}\right) \circ\left(\tilde{\mathbf{x}}_{i+1}^{(k)}-\mathbf{s}_{i+\frac{1}{2}}^{(k)}\right)+\frac{4 \mu^{(k)}}{ \beta}}\right] \label{update:x} \\
    \tau_{i+1}^{(k)} &= \arg\min_{\tau}  \left[ - \frac{\mu^{(k)}}{\beta} \log{\tau} + \frac{1}{2} \left\| \tau - \tilde{\tau}^{(k)}_{i+1} + \kappa^{(k)}_{i+\frac{1}{2}}  \right\|^{2}  \right] \nonumber \\
    &= \frac{1}{2}\left[\left(\tilde{\tau}_{i+1}^{(k)}-\kappa_{i+\frac{1}{2}}^{(k)}\right)+\sqrt{\left(\tilde{\tau}_{i+1}^{(k)}-\kappa_{i+\frac{1}{2}}^{(k)}\right) \circ\left(\tilde{\tau}_{i+1}^{(k)}-\kappa_{i+\frac{1}{2}}^{(k)}\right)+\frac{4 \mu^{(k)}}{\beta}}\right] \label{update:tau} \\
    \theta^{(k)}_{i+1} &= \tilde{\theta}^{(k)}_{i+1} - \xi^{(k)}_{i+\frac{1}{2}} - (n+1) \label{update:theta} \ 
\end{align}
The details of the half-update strategy are presented in Algorithm \ref{Alg:half-update}. We observe that this strategy can 
reduce the number of ADMM iterations on some specific datasets.

\begin{algorithm}[ht]
\SetAlgoNoLine
\caption{ABIP with half update}
\label{Alg:half-update}
Set $\mu^{0} = \beta> 0$, $\alpha, \gamma \in (0,1) $ \;
Set $\mathbf{r}^{0}_{0} = \mathbf{y}^{0}_{0} = \bzero$, $(\mathbf{x}^{0}_{0}, \tau^{0}_{0}, \mathbf{s}^{0}_{0}, \kappa^{0}_{0} ) =  ( \mathbf{e}, 1, \mathbf{e}, 1) > \bzero$, $\theta^{0}_{0} = 1$, $\xi^{0}_{0} = -n-1$ with $\mathbf{x}^{0}_{0} \circ \mathbf{s}^{0}_{0} = \frac{\mu^{0}}{\beta} \mathbf{e} $, and $\tau_{0}^{0} \kappa_{0}^{0} = \frac{\mu^{0}}{\beta}$ \;
\For{$k=0, 1, 2,\cdots$}{
    \For{$i = 0, 1, 2, \cdots$}{
        \If{the \textbf{inner} termination condition is satisfied}{
        break \;}
        Update $\Tilde{\mathbf{u}}_{i+1}^{(k)}$ by \eqref{cp-proj}  \;
        Update $ \mathbf{v}_{i+\frac{1}{2}}^{(k)} = \alpha\mathbf{u}_i^{(k)}+\textbf{v}_i^{(k)} -\alpha\Tilde{\mathbf{u}}_{i+1}^{(k)}$ \;
        Update $\mathbf{u}_{i+1}^{(k)}$ by \eqref{update:y}, \eqref{update:x}, \eqref{update:tau} and \eqref{update:theta} \;
        Update $\mathbf{v}_{i+1}^{(k)} = \mathbf{v}_{i+\frac{1}{2}}^{(k)}- \Tilde{\mathbf{u}}_{i+1}^{(k)}+\mathbf{u}_{i+1}^{(k)}$ \;
        \If{the \textbf{final} termination criterion is satisfied}{
        \Return{}\;}
}
Set $\mu^{(k+1)} = \gamma \mu^{(k)}$ \;
Set $\mathbf{r}_{0}^{(k+1)} = \bzero$, $\xi_{0}^{(k+1)} = -n-1$ and
\begin{equation*}
\left(\mathbf{y}_{0}^{(k+1)}, \mathbf{x}_{0}^{(k+1)}, \mathbf{s}_{0}^{(k+1)}, \tau_{0}^{(k+1)}, \kappa_{0}^{(k+1)}, \theta_{0}^{(k+1)}\right)=\sqrt{\gamma} \cdot\left(\mathbf{y}_{i+1}^{(k)}, \mathbf{x}_{i+1}^{(k)}, \mathbf{s}_{i+1}^{(k)}, \tau_{i+1}^{(k)}, \kappa_{i+1}^{(k)}, \theta_{i+1}^{(k)}\right)
\end{equation*}
}
\end{algorithm}

\subsection{Presolve and preconditioning}

The presolve procedure plays an essential role in modern commercial LP solvers. It analyzes the LP problem before submitting it to the optimization algorithm. With the procedure, we can detect and remove empty rows, empty columns, singleton rows and fixed variables, together with removing all the linearly dependent rows. We integrate PaPILO~\citep{Papilo}, an open-source presolve package for linear optimization, into ABIP in the enhanced implementation.

After Presolve,  we transform \eqref{eq:lp} to the following reduced LP problem:
\begin{equation*}
    \begin{split}
        \min & \quad  \bar{\bc}^{T} \bar{\bx} \\
        \text{s.t.} & \quad \bar{\bA} \bar{\bx} = \bar{\bb}, \\
                    &\quad \bar{\bx} \geq \bzero. \\
    \end{split}
\end{equation*}
Then, similar to PDLP \citep{PDLP-2021}, we perform diagonal preconditioning, which is a popular heuristic to improve the convergence of optimization algorithms. We rescale the constraint matrix $\bar{\bA} \in \mathbb{R}^{m\times n}$ to $\Tilde{\bA} = \bD_{1}^{-1}\bar{\bA}\bD_{2}^{-1}$ with positive diagonal matrices $\bD_{1}$ and $\bD_{2}$. Such preconditioning creates a new LP instance that replaces $\bar{\bA}, \bar{\bc}, \bar{\bx}$ and $\bar{\bb}$ with $\tilde{\bA}$, $ \tilde{\bc} = \bD_{2}^{-1}\bar{\bc}$, $\tilde{\bx} = \bD_{2} \bar{\bx} $ and $\tilde{\bb} = \bD_{1}^{-1}\bar{\bb}$, respectively. In particular,
we consider the following two rescaling methods.

\paragraph{Pock-Chambolle rescaling \citep{PockChambolle-2011}.} This rescaling method is parameterized by $\alpha > 0$. The diagonal matrices are defined by $(\bD_{1})_{jj} = \sqrt{\| \bA_{j, \cdot} \|_{2-\alpha}}$ for $j = 1, 2, \cdots, m $, and $(\bD_{2})_{ii} = \sqrt{\| \bA_{\cdot, i} \|_{\alpha}}$ for $i = 1, 2, \cdots, n$.
    
\paragraph{Ruiz rescaling \citep{Ruiz-2001}.} In each iteration of Ruiz rescaling, the diagonal matrices are defined by $(\bD_{1})_{jj} = \sqrt{\| \bA_{j, \cdot} \|_{\infty}}$ for $j = 1, 2, \cdots, m $, and $(\bD_{2})_{ii} = \sqrt{\| \bA_{\cdot, i} \|_{\infty}}$ for $i = 1, 2, \cdots, n$. When this rescaling method is applied iteratively, the infinity norm of each row and each column will finally converge to $1$. 

In the implementation, we first take a few iterations of Ruiz rescaling and
then apply the Pock-Chambolle rescaling with $\alpha = 1$. Then, we solve the new LP instance with ABIP, and transform its optimal solution with stored $\bD_1$ and $\bD_2$ to obtain the true one.

\subsection{Ignoring dual feasibility for null objective problems}
\label{subsection:null-obj}
In this subsection, we consider the following null objective
linear optimization problem:
\begin{eqnarray*}
  \min_{\bx} & 0  \\
  \text{s.t.} & \bA \bx = \bb,  \\
  & \bx \geq  \mathbf{0}, 
\end{eqnarray*}
which is essentially a feasibility problem. This class of problem has wide applications in real life, such as PageRank. It is used by Google to order search engine results, and aims to find the maximal right eigenvector of a stochastic matrix. \cite{nesterov-subgradient} shows that, it can be reformulated as a linear optimization problem with the null objective,  
\begin{equation*}
    \begin{split}
        \min & \quad  0\\
        \text{s.t.} & \quad \boldsymbol{S} \bx \leq \bx, \\
                    & \quad \boldsymbol{1}^{T} \bx = 1, \\
                    & \quad \bx \geq \bzero. \\
    \end{split}
\end{equation*} 

When ABIP is invoked to solve the above feasibility problem,
the dual is homogeneous and inherently admits a trivial feasible solution $(\by, \bs) = ( \mathbf{0},
 \mathbf{0})$. Therefore, we can ignore dual feasibility checks during iterations. 
Although the above observation cannot be used to identify the support of the optimal $\bs^{\ast}$,
it implies that we do not need to check dual infeasibility when the algorithm
iterates to obtain a near-feasible primal-dual pair in ABIP.

\subsection{Integrating the enhancements}

\begin{figure}[ht]
\centering
\includegraphics[scale=0.5]{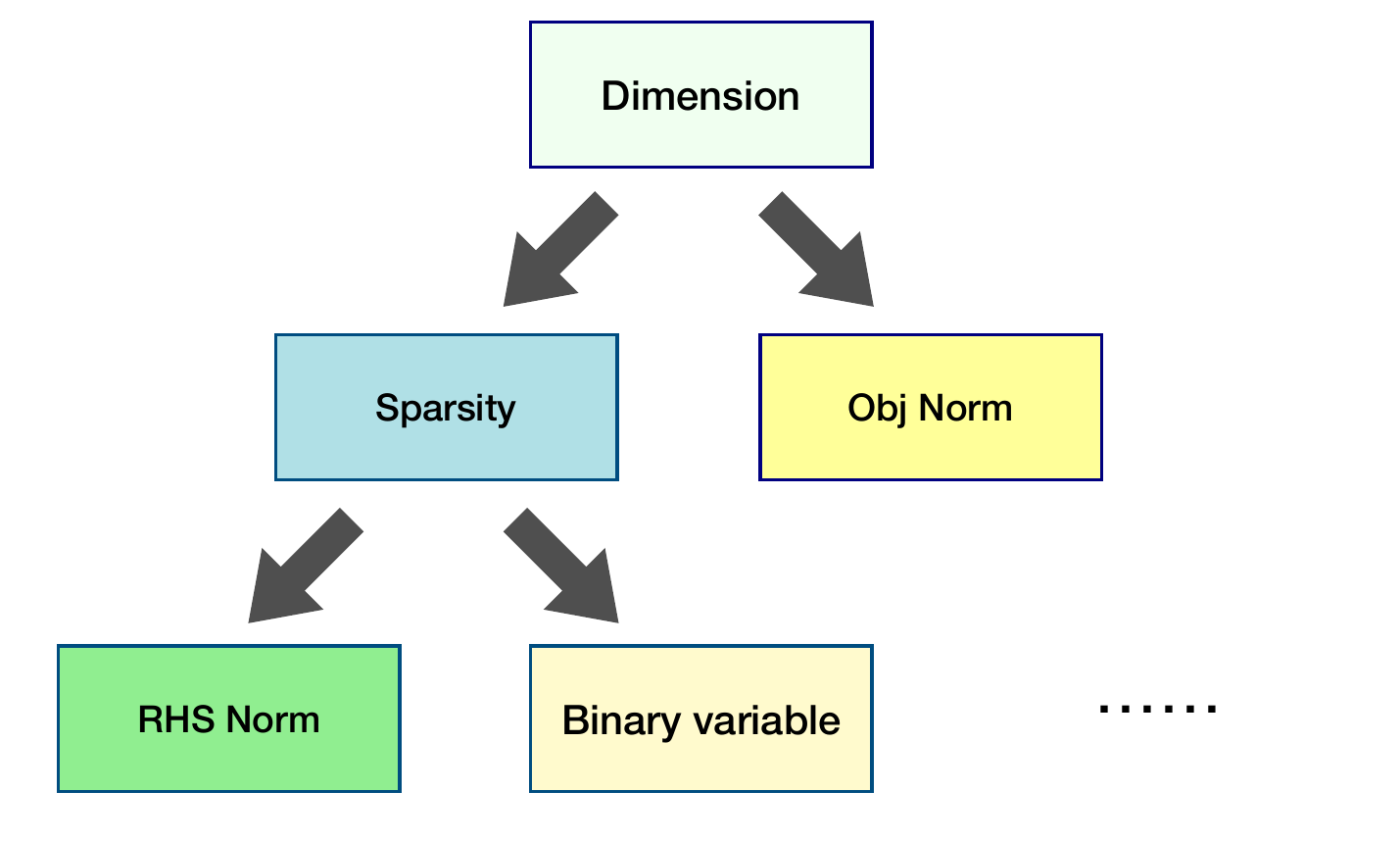}
\caption{Decision tree for ABIP}
\label{fig:decision-tree}
\end{figure}

As mentioned above, we propose various new strategies to
improve the empirical performance of ABIP.
Each strategy has some parameters to choose from, 
and thus one critical problem is how to determine the combination of parameters in
these strategies to obtain a better empirical performance. To address this challenge, a
decision tree is embedded in ABIP to help identify the rules for choosing a proper combination of
strategies. More specifically, for each LP instance, we collect the sparsity and
dimension-related statistics of the problem as features, as shown in Figure \ref{fig:decision-tree}, and record the runtime of ABIP under different combinations of strategies. Then, the combination of strategies with the least runtime is chosen to be the label of the LP instance.

To avoid overfitting, we only consider $3$ different combinations of strategies and construct the training set based on different LP datasets including \texttt{MIPLIB
2017} \citep{gleixner2021miplib} and \texttt{Netlib LP} \citep{koch2003final}. In addition, after the decision tree is obtained from the
training procedure, we prune the tree till the rules left are general enough to be coded by the solver. The computational experiments in the following section suggest that the decision tree is often able to capture features for different LP problems, and the derived decision rules prove effective in practice.

\subsection{Predictor-corrector method}

In practice, we observe that the sublinear convergence rate of ABIP often results in a ``tailing-off'' effect, preventing the algorithm from achieving satisfactory tolerance levels. This is particularly noticeable in medium-sized problems, where $10$ to $100$ ADMM iterations can be as computationally expensive as a full Newton step in traditional IPM.
To partially address this, we incorporate a primal-dual IPM with Mehrotra's corrector step into ABIP. Specifically, with a pre-specified tolerance level, if ABIP observes the significant tailing-off effect and the tolerance level is not reached, it will perform several interior point steps to improve solution accuracy. Since each interior point step involves a matrix decomposition, we only consecutively use it at most five times and employ Mehrotra's corrector step to reuse the decomposition.

\section{Numerical experiments}\label{Sec:numerical}
In this section, we present a comprehensive numerical study on the enhanced version of ABIP, referred to as ABIP+. ABIP+ is now implemented in C language with a MATLAB interface. As mentioned before, it uses two different methods to solve the linear system. The direct method uses a sparse permuted $\textsf{LDL}^{\top}$ factorization provided by the SuiteSparse package \citep{suitesparse}, while the indirect ABIP+ uses a preconditioned conjugate gradient method. 

For LP, we comprehensively evaluate the strategies of ABIP+ using the Netlib LP benchmark dataset. Then, we test the performance of ABIP+ on $240$ instances from MIP2017, $48$ instances from the Mittlemann LP barrier dataset, and randomly generated PagaRank instances. We also compare ABIP+ with PDLP \citep{PDLP-2021}, a state-of-the-art first-order LP solver implemented in C++\footnote{For more details, see \url{https://developers.google.com/optimization/lp/pdlp\_math}.}, and COPT, a general-purpose commercial solver with leading performance on LP. 

To further demonstrate the performance of ABIP+ on conic optimization, we first randomly generate several LASSO problems and choose $6$ large SVM problems from the LIBSVM dataset~\citep{chang2011libsvm}. For each problem, we maintain its two equivalent formulations. The first is the SOCP formulation. The second is the QP formulation with the quadratic objective function and conic constraints. We provide the details in Appendix \ref{sec:linsys}. Although both LASSO and SVM problems can be solved by customized algorithms developed in machine learning, it is a standard practice \citep{stellato2020osqp, SCS-2016} to use them to examine the scalability and efficiency of conic optimization and QP solvers. For ABIP+, despite that we only discuss its extension to linear conic optimization before, it can be further developed to solve the following quadratic conic optimization (QCP) and thus handle the QP formulation of LASSO and SVM problems:
\begin{equation}
\label{prob:qcp}
	\begin{aligned}
		\min & \quad \frac{1}{2}\mathbf{x}^T\mathbf{P}\mathbf{x} + \mathbf{c}^{T}\mathbf{x} \\
	\quad	\text{s.t.} &\quad \mathbf{A} \mathbf{x} = \mathbf{b}, \\
		&\quad \bx \in \mathcal{K}.
	\end{aligned} %
\end{equation}
To be succinct, ABIP+ can maintain a similar update rule as \eqref{cp-proj}, \eqref{cp-prox}, and \eqref{cp-dual} by considering the monotone linear complementarity problem and homogeneous embedding \citep{andersen1999homogeneous,o2021operator}, and we omit the details in this paper. We use the state-of-the-art open-source conic optimization solver SCS~\citep{SCS-2016} and
the general-purpose commercial solver  GUROBI~\citep{gurobi} as the benchmarks. We do not compare with OSQP~\citep{stellato2020osqp}, since it uses a weaker termination criterion and has an inferior performance to SCS as shown in \cite{o2021operator}. Finally, we assess the performance of ABIP+ on standard SOCP benchmark CBLIB~\citep{friberg2016cblib}, and compare with SCS and another recently popular first-order solver COSMO~\citep{garstka2021cosmo}.

Let $\epsilon > 0 $ be the tolerance level of the solution accuracy. We use \eqref{eq:three-all-small} as the global convergence criterion. For linear optimization, we compute $ \epsilon_{\texttt{pres}}^{\prime}, \epsilon_{\texttt{dres}}^{\prime}$, and $\epsilon_{\texttt{dgap}}^{\prime} $ by \eqref{outer:criteria}.
While for general conic optimization, they are computed by 
\begin{align*}
    \epsilon_{\texttt{pres}}^{\prime} &= \left\|\bA\frac{\bx}{\tau} - \bb\right\|_\infty / \left(1 + \max\left(\left\|\bA\frac{\bx}{\tau}\right\|_\infty, \|\bb\|_\infty\right)\right), \\
    \epsilon_{\texttt{dres}}^{\prime} &= \left\|\bA^T\frac{\by}{\tau} + \frac{\bs}{\tau} - \bc\right\|_\infty / \left(1 + \|\bc\|_\infty \right), \\
    \epsilon_{\texttt{dgap}}^{\prime} &= 
 \left|\bc^T\frac{\bx}{\tau} - \bb^T\frac{\by}{\tau}\right| /  \left(1 + \max\left(\left|\bc^T\frac{\bx}{\tau}\right|, \left|\bb^T\frac{\by}{\tau}\right|\right)\right).
\end{align*}
Typically, we set $\epsilon = 10^{-4}$ or $10^{-6}$. For inner loop in conic optimization, we use the following termination criteria: $||\bQ\bu_i^{(k)}-\bv_i^{(k)}|| \leq (\mu^{(k)})^\alpha(1+||(\bu_i^{(k)},\bv_i^{(k)})||)$ with $\alpha \in [0.25,2]$, which ensures the satisfactory convergence behavior of ABIP+.

When testing the Netlib LP benchmark dataset, we limit the number of ADMM iterations to $10^{6}$. For other LP benchmark datasets, we set the runtime limit to be $3600$s. For conic optimization, we set the runtime limit to be $2000$s for LASSO, $3000$s for SVM, and $100$s for SOCP in CBLIB. Finally, we remark that we do not apply any of the novel strategies developed in Section \ref{Sec:new-strategy} in the general conic optimization since the acceleration performance is marginal in the experiments.

Unless otherwise mentioned, we use the shifted geometric mean (SGM) of runtimes to evaluate the algorithm's performance over a dataset. Specifically, we define the SGM by 
\begin{equation*}
    \texttt{SGM} = \left( \prod_{i=1}^{n} \left(t_{i} + \Delta \right) \right)^{1/n} - \Delta,
\end{equation*}
where $t_{i}$ is the runtime for the $i$-th instance, and $\Delta = 10$ is a time-shift value to alleviate the effect of runtimes that are almost $0$. When the algorithm fails to solve a certain instance, we set $t_{i} = 15000$. We normalize all the SGMs by setting the smallest SGM to be $1$.

\subsection{Numerical results on Netlib LP benchmark dataset}
In this subsection, we showcase the improvements of ABIP+ over the original ABIP on the Netlib LP benchmark dataset. For the restart strategy, we set the restart threshold $M_T=10^5$ and the fixed restart period $F = 10^3$. For the rescale strategy, we do Pock-Chambolle rescaling first and then Ruiz rescaling $10$ times. For the hybrid $\mu$ strategy, we use the aggressive strategy when the barrier parameter $\mu^{(k)} < 10^3 \cdot \epsilon$. Otherwise, we use the LOQO method.

Table \ref{table:netlib} shows the performance of ABIP with different strategies. The first row gives the results of the original ABIP, and the second row shows the performance of ABIP with restart strategy. Similarly, the third and the fourth row correspond to ABIP with both restart and rescale strategies and ABIP with the whole three strategies, respectively. The second column gives the number of instances each algorithm solved within $10^6$ ADMM iterations. Moreover, the last three columns give the average number of IPM iterations, the average number of ADMM iterations, and the average computational time over the solved instances, respectively. 
It can be seen that each added strategy improves the performance of ABIP. Among $105$ LP instances selected from the Netlib dataset, ABIP solves only $65$ instances to $10^{-6}$ relative accuracy given a limit of $10^{6}$ ADMM iterations per problem. On the contrary, ABIP+ solves $86$ instance under the same setting. Moreover, ABIP+ reduces more than $70\%$ ADMM iterations and more than $80\%$ runtime on average. Therefore, the acceleration performance of these strategies is significant.
\begin{table}[htbp]
        \centering
        \begin{tabular}{l|c|c|c|c}
            \hline
            Method     &   \# Solved   & Avg. \# IPM & Avg. \# ADMM   &  Avg.Time (s)      \\
            \hline
            ABIP       &     $65$           &   $74$    &   $265418$     &   $87.07$                    \\
            + restart   &   68 & 74 &  88257 & 23.63  \\
            + rescale   &  84 &  72 & 77925 &  20.44 \\
            + hybrid $\mu$ (=ABIP+) & $\boldsymbol{86}$ & $\boldsymbol{22}$ & $ \boldsymbol{73738}$ & $ \boldsymbol{14.97} $ \\
            \hline
        \end{tabular}
        \caption{Performance of ABIP with different strategies on Netlib, $\epsilon = 10^{-6}$.}
        \label{table:netlib}
    \end{table}

\subsection{Numerical results on MIP2017 LP benchmark dataset}\label{MIP2017}
In this subsection, we evaluate the performance of ABIP+, PDLP and COPT, on the MIP2017 LP benchmark dataset.
For a fair comparison, all tested instances for PDLP, ABIP, and ABIP+ are presolved by PaPILO \citep{Papilo}. We also include different tolerance levels to terminate the first-order LP solvers. Specifically, we let $\epsilon = 10^{-4}$ and $10^{-6}$ for PDLP and ABIP+. Since ABIP+ inherits the framework of IPM, we also test the effectiveness of the predictor-corrector method. Specifically, after the ADMM iterates reach the tolerance level $\epsilon=10^{-3}$, we invoke the predictor-corrector method to further reach $\epsilon=10^{-6}$. These results are indicated with a $\dagger$ symbol presented in Table \ref{table:mip2017} and Table \ref{table:mip2017-instance}.

On one hand, Table \ref{table:mip2017} shows that the first-order algorithms generally lag behind the IPM implemented in COPT. On the other hand, it demonstrates that ABIP+ significantly improves the performance of original ABIP. It is also comparable to PDLP on the MIP2017 dataset when $\epsilon=10^{-6}$ and the  predictor-corrector method is used. 
While PDLP may slightly outperform ABIP+ in terms of overall performance, 
Table \ref{table:mip2017-instance} gives specific instances where ABIP+ has a clear advantage over PDLP in terms of computational time. 
This suggests that ABIP+ and PDLP may each be better suited to different types of problems, which may require independent investigation.
\begin{table}[htbp]
    \centering
    \begin{tabular}{l|c|c|c}
    \hline
    Method  &  $\epsilon$ & \# Solved  &  SGM  \\
    \hline
    COPT    & $10^{-6}$  & $\boldsymbol{240}$   &   $\boldsymbol{1}$      \\ 
    \hline
    ABIP & $10^{-4}$ &  192   & 34.8        \\
    \hline
    \multirow{2}{*}{ABIP+} & $10^{-4}$ &${220}$   &  ${21.3}$       \\
             & $10^{-6}$ &${213}$   &  ${35.4}$       \\
             & $10^{-6^\dagger}$ &${220}$   &  ${14}$       \\
    \hline
    \multirow{2}{*}{PDLP}    & $10^{-4}$ & $226$ &  $6.7$         \\
            & $10^{-6}$ & $221$ &  $12.9$         \\
    \hline
    \end{tabular}
    \caption{Numerical results on MIP2017. The symbol $^\dagger$ means that we utilize the predictor-corrector method.}
    \label{table:mip2017}
\end{table}

\begin{table}[htbp]
    \centering
    \begin{tabular}{l|c|c|c|c}
    \hline
    Instance     &  ABIP+ & ABIP+$^\dagger$ &  PDLP  & COPT   \\
    \hline
    \texttt{app1-2}  &  49.65  & 8.07 & 118.83   &  \textbf{1.18}  \\
    \texttt{buildingenergy} & 508.34 & 74.98 &  998.21   &  \textbf{18.09}    \\
    \texttt{map16715-04} & 638.13  & 17.57 &  1241.83  & \textbf{3.08}    \\
    \texttt{neos-4647030-tutaki}  & 403.03  &266.87&  f  &    \textbf{4.71}   \\
    \texttt{neos-4763324-toguru}  & 1336.10  &997.77&  f  &  \textbf{6.06} \\
    \texttt{unitcal\_7}    &   27.12   & 32.39 & 282.34  & \textbf{3.61} \\
    \hline
    \end{tabular}
    \caption{The runtime (in sec.) of ABIP+, ABIP+$^\dagger$, PDLP and COPT on selected instances from MIP2017. For ABIP+ and PDLP, we set $\epsilon = 10^{-4}$. For ABIP+$^\dagger$ and COPT, we set $\epsilon = 10^{-6}$. The symbol $^\dagger$ means that we utilize the predictor-corrector method. Note that ``f'' in the table means that the algorithm fails to solve the instance within the runtime limit. }
    \label{table:mip2017-instance}
\end{table}

\subsection{Numerical results on Mittlemann LP barrier benchmark dataset}
To further test the robustness of different solvers studied in the last subsection, we compare their performance on the Mittlemann LP barrier dataset. This dataset is commonly used for benchmarking IPM-based LP solvers. 
As before, we include the results of ABIP+ utilizing the predictor-corrector method, marked by the symbol $\dagger$, at a tolerance level of $\epsilon=10^{-6}$. Similarly, the barrier IPM in COPT is also terminated at the same accuracy.
Table \ref{table:mittleman-barrier} shows that, although our ABIP+ is slightly inferior to PDLP in terms of computational time, it solves more instances at $\epsilon=10^{-6}$. Meanwhile, ABIP+ with the predictor-corrector method also shows its effectiveness in this dataset. Table \ref{table:mittleman-barrier-instance} lists some instances that ABIP+ solves faster than PDLP. 
\begin{table}[htbp]
    \centering
    \begin{tabular}{l|c|c|c}
    \hline
    Method  & $\epsilon$ & \# Solved  &  SGM  \\
    \hline
    COPT    & $10^{-6}$  & $\boldsymbol{50}$   &   $\boldsymbol{1}$      \\ 
    \hline
    ABIP & $10^{-4}$ &  22  & 57.40         \\
    \hline
    \multirow{2}{*}{ABIP+} & $10^{-4}$ &${29}$   &  ${27.22}$       \\
             & $10^{-6}$ &${22}$   &  ${49.32}$       \\
             & $10^{-6^\dagger}$ &${32}$   &  ${17.22}$       \\
    \hline
    \multirow{2}{*}{PDLP}    & $10^{-4}$ & $39$ &  $9.57$         \\
            & $10^{-6}$ & $19$ &  $34.54$         \\
    \hline
    \end{tabular}
    \caption{Numerical results on Mittlemann LP barrier dataset.}
    \label{table:mittleman-barrier}
\end{table}
\begin{table}[htbp]
    \centering
    \begin{tabular}{l|c|c|c}
    \hline
    Instance     &  ABIP+  &  PDLP  & COPT   \\
    \hline
    \texttt{cont1}  &  8.93  &  49.17  &  \textbf{3.35}  \\
     \texttt{irish-e}  &  31.61    &    536.14  &   \textbf{24.81}   \\
     \texttt{ns1688926}  &   1200.622  &  f  &  \textbf{19.06}    \\
     \texttt{rmine15}  &   713.36   &   2782.00 &  \textbf{87.99} \\
\hline
    \end{tabular}
    \caption{The runtime (in sec.) of ABIP+, PDLP and COPT on selected instances from the Mittlemann LP barrier dataset. For ABIP+ and PDLP, we set $\epsilon=10^{-4}$.}
    \label{table:mittleman-barrier-instance}
\end{table}

\subsection{Numerical results on PageRank instances}
\label{subsec:numerical-pagerank}
In this subsection, we apply ABIP+ to solve the PageRank problem via its LP formulation, as introduced in Section \ref{subsection:null-obj}. Due to a special sparse pattern in the coefficient matrix, the traditional IPM solvers, such as GUROBI, fail to solve most of these instances, as shown in \cite{PDLP-2021}. The reason behind is that, this special sparse pattern takes much time in the ordering stage and consumes a massive amount of memory (typically above 30GBytes). However, since ABIP+ can utilize the conjugate gradient method instead of direct factorization for these instances, together with the strategies introduced in Section \ref{Sec:new-strategy}, our ABIP+ manages to solve these PageRank instances efficiently. 

The PageRank instances are first benchmarked by PDLP \citep{PDLP-2021}. Following the instructions in \cite{PDLP-2021}, we generate 115 PageRank instances from several sparse matrix datasets, including DIMACS10, Gleich, Newman, and SNAP. These matrices can be found from the SuiteSparse matrix collection\footnote{For details, see \url{https://sparse.tamu.edu/about}.}.
Besides a comparison over runtime, we also compare the number of matrix-vector products used in these methods. Table \ref{table:pagerank-typical-instance} provides some specific instances where ABIP+ can be significantly faster than PDLP.

\begin{table}[htbp]
    \centering
    \begin{tabular}{l|c|c|c}
    \hline
    Method  &  \# Solved  &  SGM & Scaled Mean Matrix-Vector Products \\
    \hline
    ABIP+ &  $\boldsymbol{107}$   &  $\boldsymbol{1}$  & $\boldsymbol{1528.72}$ \\
    PDLP    &   $101$ &  $3.91$ & $6510.48$      \\
    \hline
    \end{tabular}
    \caption{ABIP+ v.s. PDLP on randomly generated PageRank dataset, $\epsilon = 10^{-4}$.}
    \label{table:pagerank-dataset}
\end{table}	
\begin{table}[htbp]
    \centering
    \begin{tabular}{l|c|c|c}
    \hline
    Instance    & \# nodes  &  PDLP  &  ABIP+  \\
    \hline
    \texttt{usroads}  &  129164  &  51.66   &   $ \boldsymbol{7.16} $ \\
    \texttt{vsp\_bcsstk30\_500sep\_10in\_1Kout}  & 58348  & 252.86  &  $\boldsymbol{10.72}$  \\
    \texttt{web-BerkStan}   &  685230  & 3497.32  &     $\boldsymbol{552.44}$     \\
    \hline
    \end{tabular}
    \caption{The runtime (in sec.) of selected instances from randomly generated PageRank dataset, $\epsilon = 10^{-4}$.}
    \label{table:pagerank-typical-instance}
\end{table}

Another interesting class of PageRank instances where ABIP+ has a clear advantage is the staircase PageRank instances. When one generates the PageRank instance by the code provided in \cite{PDLP-2021}, the coefficient matrix will exhibit a staircase form if the number of nodes is set to equal to the number of edges, see Figure \ref{fig:staircase}. In Table \ref{table:pagerank-staircase}, we list several such instances and compare ABIP+ with PDLP. It is notable that, as the matrix size increases, the runtime gap between PDLP and ABIP+ becomes more prominent. 

\begin{figure}
    \centering
    \begin{equation*}
       \begin{pmatrix}
        -1 & 0.198 & 0 & 0 & 0 & 0 & 0 & 0 & 0 & 0 \\
        0.99 & -1 & 0.495 & 0.99 & 0.495 & 0.495 & 0 & 0 & 0 & 0 \\
        0 & 0.198 & -1 & 0 & 0 & 0 & 0.495 & 0 & 0 & 0 \\
        0 & 0.198 & 0 & -1 & 0 & 0 & 0 & 0 & 0 & 0 \\
        0 & 0.198 & 0 & 0 & -1 & 0 & 0 & 0.99 & 0 & 0 \\
        0 & 0.198 & 0 & 0 & 0 & -1 & 0 & 0 & 0.99 & 0 \\
        0 & 0 & 0.495 & 0 & 0 & 0 & -1 & 0 & 0 & 0.99 \\
        0 & 0 & 0 & 0 & 0.495 & 0 & 0 & -1 & 0 & 0 \\
        0 & 0 & 0 & 0 & 0 & 0.495 & 0 & 0 & -1 & 0 \\
        0 & 0 & 0 & 0 & 0 & 0 & 0.495 & 0 & 0 & -1
        \end{pmatrix}
    \end{equation*}
    \caption{An Illustration of coefficient matrix of the staircase PageRank instance with $10$ nodes.}
    \label{fig:staircase}
\end{figure}

\begin{table}[htbp]
    \centering
    \begin{tabular}{c|c|c}
    \hline
    \# nodes &  PDLP           &  ABIP+    \\
     \hline
      $10^{4}$       &   $8.60$           &   $ \boldsymbol{0.93} $      \\
      $10^{5}$       &  $135.67$          &   $ \boldsymbol{10.36}$              \\
      $10^{6}$       &  $2248.40$       &  $\boldsymbol{60.32}$             \\
      \hline
    \end{tabular}
    \caption{Staircase PageRank instances, $\epsilon = 10^{-6}$.}
    \label{table:pagerank-staircase}
\end{table}

\subsection{Numerical results on LASSO instances}
From this subsection, we evaluate the performance of ABIP+ on conic optimization. In particular,
we apply ABIP+ to solve LASSO, SVM and standard SOCP problems. In this section, we focus on LASSO problem, which is a well-known regression model obtained by adding an $\ell_{1}$-norm regularization term in the objective to promote sparsity. It can be formulated as
\[
\min_\bx \ \|\bA\bx-\bb\|_{2}^{2} + \lambda\|\bx\|_{1},
\]
where $\bA\in\mathbb{R}^{m\times n}$ is the data matrix and $\lambda$ is the $\ell_1$-penalty weight. Here, $m$ is the number of samples and $n$ is the number of features. As mentioned before, we solve the problem via its SOCP and QCP formulations. For ABIP+, we also utilize the customized linear system solver by exploiting the intrinsic structure of the problem.

We randomly generate 9 LASSO problems, which span a wide range of ratios between $m$ and $n$. We choose $m \in \{1000,2000,5000\}, n\in\{5000,10000,15000\}$. In all the experiments, we let $\lambda=\|\bA^T \bb\|_\infty/5$, the tolerance level $ \epsilon = 10^{-3}$, and set the runtime limit to be 2000 seconds. 
We compare the performance of ABIP+ with SCS and GUROBI, and provide the results in Table~\ref{table:Lasso}. It shows that ABIP+ with SOCP formulation outperforms SCS and GUROBI with both formulations over all the instances, suggesting the empirical advantage of ABIP+ over other solvers.
Moreover, it is interesting to observe that GUROBI almost has the lowest iteration number among all the three solvers. However, its overall runtime performance is still inferior to ABIP+ due to the expensive computation cost when repetitively solving the Newton equation.

\begin{table}[htbp]
  \centering
    \begin{tabular}{cc|cc|cc|cc|cc|cc|cc}
    \hline
    \multicolumn{2}{c|}{} & \multicolumn{2}{c|}{ABIP+} & \multicolumn{2}{c|}{ABIP+} & \multicolumn{2}{c|}{SCS} & \multicolumn{2}{c|}{SCS} & \multicolumn{2}{c|}{GUROBI} & \multicolumn{2}{c}{GUROBI} \\
    \multicolumn{2}{c|}{} & \multicolumn{2}{c|}{SOCP} & \multicolumn{2}{c|}{QCP} & \multicolumn{2}{c|}{SOCP} & \multicolumn{2}{c|}{QCP} & \multicolumn{2}{c|}{SOCP} & \multicolumn{2}{c}{QCP} \\
    \hline
    $m$     & $n$     & Time  & Iter  & Time  & Iter  & Time  & Iter  & Time  & Iter  & Time  & Iter  & Time  & Iter \\
    \hline
    1000  & 5000  & $\boldsymbol{1.11}$  & 67    & 6.27  & 387   & 4.31  & 675   & 6.50  & 250   & 1.19  & $\boldsymbol{8}$     & 1.28  & $\boldsymbol{8}$ \\
    1000  & 10000 & $\boldsymbol{1.95}$  & 61    & 21.73 & 820   & 16.09 & 1275  & 16.95 & 275   & 2.76  & $\boldsymbol{9}$     & 2.69  & $\boldsymbol{9}$ \\
    1000  & 15000 & $\boldsymbol{3.62}$  & 60    & 47.18 & 1179  & 22.33 & 1600  & 23.09 & 225   & 4.23  & $\boldsymbol{10}$    & 4.02  & $\boldsymbol{10}$ \\
    2000  & 5000  & $\boldsymbol{1.92}$  & 35    & 11.68 & 335   & 12.79 & 525   & 17.80 & 175   & 3.08  & 9     & 2.87  & $\boldsymbol{7}$ \\
    2000  & 10000 & $\boldsymbol{4.17}$  & 35    & 33.02 & 606   & 29.78 & 925   & 32.64 & 200   & 6.77  & 11    & 5.75  & $\boldsymbol{8}$ \\
    2000  & 15000 & $\boldsymbol{6.52}$  & 39    & 72.45 & 898   & 43.23 & 1125  & 48.20 & 200   & 9.53  & $\boldsymbol{9}$     & 8.86  & $\boldsymbol{9}$ \\
    5000  & 5000  & $\boldsymbol{21.15}$ & 161   & 31.50 & 276   & 134.17 & 1250  & 140.50 & 150   & \multicolumn{2}{c|}{timeout}   & 10.95 & $\boldsymbol{7}$ \\
    5000  & 10000 & $\boldsymbol{19.80}$ & 32    & 78.30 & 476   & 176.02 & 625   & 193.27 & 150   & 49.40 & 12    & 31.05 & $\boldsymbol{7}$ \\
    5000  & 15000 & $\boldsymbol{26.80}$ & 31    & 133.45 & 597   & 208.23 & 750   & 231.44 & 150   & 42.25 & 12    & 34.65 & $\boldsymbol{8}$ \\
    \hline
    \multicolumn{2}{c|}{SGM} & \multicolumn{2}{c|}{$\boldsymbol{1.00}$} & \multicolumn{2}{c|}{4.83} & \multicolumn{2}{c|}{5.59} & \multicolumn{2}{c|}{6.22} & \multicolumn{2}{c|}{3.05} & \multicolumn{2}{c}{1.14} \\
    \hline
    \end{tabular}%
  	\caption{Comparison of runtime (in sec.) and iteration numbers on the LASSO problems.}
	\label{table:Lasso}%
\end{table}%

\subsection{Numerical results on SVM instances}
In this subsection, we compare the performance of ABIP+ with SCS and GUROBI on SVM problem. 
Specifically, we consider the binary classification problem and let $\{\bx_i,\by_i\}_{1\le i\le m}$ be the training data, where $\bx_i\in\mathbb{R}^n$ is the feature vector and $y_i \in \{1, -1\}$ is the corresponding label. SVM aims to classify the two groups by a linear model, and solves the following quadratic problem:
\begin{equation*}
\begin{aligned}
	\min_{(\bw,b)\in\mathbb{R}^{n+1}}\quad & \frac{\lambda}{2}\bw^T\bw + \frac{1}{m}\sum_{i=1}^m \xi_i \\
	\text{s.t.}\quad & y_i(\bx_i^T\bw + b) \geq 1 - \xi_i, \quad \forall i = 1,2,...,m, \\
	& \xi_i \geq 0, \quad \forall i = 1,2,...,m.
\end{aligned}
\end{equation*}

We choose $6$ large SVM instances from LIBSVM~\citep{chang2011libsvm}, which cover both scenarios of high dimensionality and large sample size. In all the experiments, we choose $\lambda=10^{-3}$, let $\epsilon = 10^{-3}$, and set the runtime limit to be 3000s. We also solve the problem via its SOCP and QCP formulations, and present the results in Table~\ref{tab:svm_socp_time}. It demonstrates that ABIP+ has the most robust performance among all the compared solvers. It successfully solves all the problems in both QCP and SOCP formulations, while SCS and GUROBI are unable to solve all instances in each formulation. Specifically, GUROBI is unable to solve the \texttt{real-sim} instance in both formulations. SCS exhibits inferior performance on large-scale problems, failing on \texttt{news20} and \texttt{real-sim} instances, and consuming significant runtime for the \texttt{rcv1\_train} instance in the SOCP formulation. Furthermore, it fails on \texttt{real-sim} and \texttt{skin\_nonskin} instances and experiences a large runtime for the \texttt{news20} instance in the QCP formulation.

\begin{table}[htbp]
  \centering
    \begin{tabular}{c|c|c|c|c|c|c|c|c}
    \hline
    Dataset & $m$ & $n$ & \makecell[c]{ABIP+ \\ SOCP} & \makecell[c]{ABIP+ \\ QCP} & \makecell[c]{SCS \\ SOCP} & \makecell[c]{SCS \\ QCP} & \makecell[c]{GUROBI \\ SOCP} & \makecell[c]{GUROBI \\ QCP} \\
    \hline
    \texttt{s} & 581012 & 54    & $\boldsymbol{31.85}$ & 270.92 & 134.98 & 67.97 & 1513.57 & 1711.41 \\
    \texttt{ijcnn1} & 49990 & 22    & $\boldsymbol{2.76}$  & 4.00     & 6.45  & 2.98 & 13.4  & 13.79 \\
    \texttt{news20} & 19996 & 1355191 & 449.23 & 209.35 & timeout & 2877.66 & $\boldsymbol{138.19}$ & 147.7 \\
    \texttt{rcv1\_train} & 20242 & 44504 & 225.42 & 168.41 & 2352.65 & 962.9 & $\boldsymbol{91.35}$ & 138.97 \\
    \texttt{real-sim} & 72309 & 20958 & 353.85 & $\boldsymbol{232.54}$ & timeout & timeout & timeout & timeout \\
    \texttt{skin\_nonskin} & 245057 & 3     & $\boldsymbol{9.27}$  & 67.72 & 132.22 & timeout & 91.48 & 91.71 \\
    \hline
    \end{tabular}%
\caption{Comparison of runtime (in sec.) on the selected SVM problems.}
\label{tab:svm_socp_time}
\end{table}%

\subsection{Numerical results on SOCP instances of CBLIB}
In this subsection, we compare the performance of ABIP+ with the other two first-order solvers SCS~\citep{SCS-2016} and COSMO~\citep{Garstka_2021} on $1405$ selected SOCP instances from the CBLIB~\citep{friberg2016cblib}. 
All the instances are presolved by COPT. We set $\epsilon = 10^{-4}$ and set the runtime limit to be $100$ seconds. When conducting experiments, we find that the first-order algorithms sometimes obtain inaccurate solutions even though their empirical stopping criteria are satisfied. For a fair comparison, when an algorithm solves an instance within the runtime limit, we count the instance as a solved instance only if its relative error to the optimal value returned by MOSEK is less than a pre-specified threshold $\epsilon_{\text{tol}}$. The relative error is defined by $(|f - f^*|) / \max\{|f^*|,1\}$, where $f^*$ is the optimal value returned by MOSEK, and we set $\epsilon_{\text{tol}}=0.01$. The detailed results are presented in Table~\ref{tab:cblib}. It can be seen that MOSEK and COPT achieve better performance than the first-order solvers, since second-order methods have the advantage of high accuracy and efficiency when solving small-scale problems. When comparing between the first-order solvers, SCS has a better performance than ABIP+ on runtime, but only solves 3 more instances. Hence, our ABIP+ has a comparable performance with SCS. Meanwhile, we find that both SCS and ABIP+ outperform COSMO in terms of solving more instances and requiring less runtime.

\begin{table}[htbp]
    \centering
      \begin{tabular}{c|ccccc}
        \hline
        & MOSEK & COPT & SCS & ABIP+ & COSMO \\
        \hline
      \# Solved & $\boldsymbol{1405}$ & 1387 & 1337 & 1334 & 1308 \\
      SGM   & $\boldsymbol{1.00}$ & 1.13 & 3.26 & 3.52 & 4.64  \\
      \hline
    \end{tabular}%
      \caption{Numerical results on CBLIB dataset.}
    \label{tab:cblib}%
\end{table}%

\section{Conclusion}
In this paper, we continue the development of the ADMM-based interior point method~\citep{ABIP-2021}. We generalize ABIP to deal with general conic constraints, theoretically justify that ABIP converges at a $\tilde{O}(1/\epsilon)$ rate for the general linear conic problems, and provide efficient implementation to further accelerate the empirical performance on some important SOCP problems. Moreover, we provide several implementation techniques that are inspired by some existing methods and substantially improve the performance of ABIP for large-scale LP. Despite its sensitivity to the condition number and inverse precision, the enhanced ABIP solver is highly competitive for certain structured challenging problems such as PageRank and some machine learning problems. We believe that ABIP+ is complementary to the state-of-the-art open-source solvers and exhibits a strong potential for future improvement.
For example, it would be interesting to extend ABIP+ to the distributed and asynchronous environment. Our code is open-source and available at \url{https://github.com/leavesgrp/ABIP}.

\bibliographystyle{unsrtnat}
\bibliography{ref.bib}

\begin{thebibliography}{50}
\providecommand{\natexlab}[1]{#1}
\providecommand{\url}[1]{\texttt{#1}}
\expandafter\ifx\csname urlstyle\endcsname\relax
  \providecommand{\doi}[1]{doi: #1}\else
  \providecommand{\doi}{doi: \begingroup \urlstyle{rm}\Url}\fi

\bibitem[Alizadeh and Goldfarb(2003)]{alizadeh2003second}
Farid Alizadeh and Donald Goldfarb.
\newblock Second-order cone programming.
\newblock \emph{Mathematical programming}, 95\penalty0 (1):\penalty0 3--51,
  2003.

\bibitem[Nesterov and Nemirovskii(1994)]{nesterov1994interior}
Yurii Nesterov and Arkadii Nemirovskii.
\newblock \emph{Interior-point polynomial algorithms in convex programming}.
\newblock SIAM, 1994.

\bibitem[Sturm(1999)]{sturm1999using}
Jos~F Sturm.
\newblock Using sedumi 1.02, a matlab toolbox for optimization over symmetric
  cones.
\newblock \emph{Optimization methods and software}, 11\penalty0 (1-4):\penalty0
  625--653, 1999.

\bibitem[Toh et~al.(1999)Toh, Todd, and T{\"u}t{\"u}nc{\"u}]{toh1999sdpt3}
Kim-Chuan Toh, Michael~J Todd, and Reha~H T{\"u}t{\"u}nc{\"u}.
\newblock Sdpt3—a matlab software package for semidefinite programming,
  version 1.3.
\newblock \emph{Optimization methods and software}, 11\penalty0 (1-4):\penalty0
  545--581, 1999.

\bibitem[{ApS, Mosek}(2019)]{aps2019mosek}
{ApS, Mosek}.
\newblock Mosek optimization toolbox for matlab, 2019.

\bibitem[{Gurobi Optimization, LLC}(2022)]{gurobi}
{Gurobi Optimization, LLC}.
\newblock {Gurobi Optimizer Reference Manual}, 2022.
\newblock URL \url{https://www.gurobi.com}.

\bibitem[Ge et~al.(2022)Ge, Huangfu, Wang, Wu, and Ye]{ge2022cardinal}
Dongdong Ge, Qi~Huangfu, Zizhuo Wang, Jian Wu, and Yinyu Ye.
\newblock Cardinal optimizer (copt) user guide.
\newblock \emph{arXiv preprint arXiv:2208.14314}, 2022.

\bibitem[Pougkakiotis and Gondzio(2022)]{pougkakiotis2022interior}
Spyridon Pougkakiotis and Jacek Gondzio.
\newblock An interior point-proximal method of multipliers for linear positive
  semi-definite programming.
\newblock \emph{Journal of Optimization Theory and Applications}, pages 1--33,
  2022.

\bibitem[Cipolla and Gondzio(2023)]{cipolla2023proximal}
Stefano Cipolla and Jacek Gondzio.
\newblock Proximal stabilized interior point methods and low-frequency-update
  preconditioning techniques.
\newblock \emph{Journal of Optimization Theory and Applications}, 197\penalty0
  (3):\penalty0 1061--1103, 2023.

\bibitem[Zhou and Toh(2004)]{zhou2004polynomiality}
Guanglu Zhou and Kim-Chuan Toh.
\newblock Polynomiality of an inexact infeasible interior point algorithm for
  semidefinite programming.
\newblock \emph{Mathematical programming}, 99:\penalty0 261--282, 2004.

\bibitem[Bellavia and Pieraccini(2004)]{bellavia2004convergence}
Stefania Bellavia and Sandra Pieraccini.
\newblock Convergence analysis of an inexact infeasible interior point method
  for semidefinite programming.
\newblock \emph{Computational Optimization and Applications}, 29:\penalty0
  289--313, 2004.

\bibitem[Lu et~al.(2006)Lu, Monteiro, and O'Neal]{lu2006iterative}
Zhaosong Lu, Renato~DC Monteiro, and Jerome~W O'Neal.
\newblock An iterative solver-based infeasible primal-dual path-following
  algorithm for convex quadratic programming.
\newblock \emph{SIAM Journal on optimization}, 17\penalty0 (1):\penalty0
  287--310, 2006.

\bibitem[Al-Jeiroudi and Gondzio(2009)]{al2009convergence}
Ghussoun Al-Jeiroudi and Jacek Gondzio.
\newblock Convergence analysis of the inexact infeasible interior-point method
  for linear optimization.
\newblock \emph{Journal of Optimization Theory and Applications}, 141:\penalty0
  231--247, 2009.

\bibitem[Zanetti and Gondzio(2023)]{zanetti2023new}
Filippo Zanetti and Jacek Gondzio.
\newblock A new stopping criterion for krylov solvers applied in interior point
  methods.
\newblock \emph{SIAM Journal on Scientific Computing}, 45\penalty0
  (2):\penalty0 A703--A728, 2023.

\bibitem[Yang et~al.(2015)Yang, Sun, and Toh]{YangST15}
Liuqin Yang, Defeng Sun, and Kim{-}Chuan Toh.
\newblock {SDPNAL} $+$ : a majorized semismooth newton-cg augmented lagrangian
  method for semidefinite programming with nonnegative constraints.
\newblock \emph{Math. Program. Comput.}, 7\penalty0 (3):\penalty0 331--366,
  2015.
\newblock \doi{10.1007/s12532-015-0082-6}.
\newblock URL \url{https://doi.org/10.1007/s12532-015-0082-6}.

\bibitem[O’Donoghue et~al.(2016)O’Donoghue, Chu, Parikh, and
  Boyd]{SCS-2016}
Brendan O’Donoghue, Eric Chu, Neal Parikh, and Stephen Boyd.
\newblock Conic optimization via operator splitting and homogeneous self-dual
  embedding.
\newblock \emph{Journal of Optimization Theory and Applications}, 169\penalty0
  (3):\penalty0 1042--1068, 2016.

\bibitem[Applegate et~al.(2021{\natexlab{a}})Applegate, Hinder, Lu, and
  Lubin]{PDHG-2021}
David Applegate, Oliver Hinder, Haihao Lu, and Miles Lubin.
\newblock Faster first-order primal-dual methods for linear programming using
  restarts and sharpness.
\newblock \emph{arXiv preprint arXiv:2105.12715}, 2021{\natexlab{a}}.

\bibitem[Boyd et~al.(2011)Boyd, Parikh, Chu, Peleato, and Eckstein]{BoydPCPE11}
Stephen~P. Boyd, Neal Parikh, Eric Chu, Borja Peleato, and Jonathan Eckstein.
\newblock Distributed optimization and statistical learning via the alternating
  direction method of multipliers.
\newblock \emph{Found. Trends Mach. Learn.}, 3\penalty0 (1):\penalty0 1--122,
  2011.
\newblock \doi{10.1561/2200000016}.
\newblock URL \url{https://doi.org/10.1561/2200000016}.

\bibitem[Ye et~al.(1994)Ye, Todd, and Mizuno]{YeToddMizuno-1994}
Yinyu Ye, Michael~J Todd, and Shinji Mizuno.
\newblock An {O} ($\sqrt{n}${L})-iteration homogeneous and self-dual linear
  programming algorithm.
\newblock \emph{Mathematics of operations research}, 19\penalty0 (1):\penalty0
  53--67, 1994.

\bibitem[Sopasakis et~al.(2019)Sopasakis, Menounou, and
  Patrinos]{sopasakis2019superscs}
Pantelis Sopasakis, Krina Menounou, and Panagiotis Patrinos.
\newblock Superscs: fast and accurate large-scale conic optimization.
\newblock In \emph{2019 18th European Control Conference (ECC)}, pages
  1500--1505. IEEE, 2019.

\bibitem[Lin et~al.(2021)Lin, Ma, Ye, and Zhang]{ABIP-2021}
Tianyi Lin, Shiqian Ma, Yinyu Ye, and Shuzhong Zhang.
\newblock An admm-based interior-point method for large-scale linear
  programming.
\newblock \emph{Optimization Methods and Software}, 36\penalty0 (2-3):\penalty0
  389--424, 2021.
\newblock \doi{10.1080/10556788.2020.1821200}.

\bibitem[Chambolle and Pock(2011)]{ChambolleP11}
Antonin Chambolle and Thomas Pock.
\newblock A first-order primal-dual algorithm for convex problems with
  applications to imaging.
\newblock \emph{J. Math. Imaging Vis.}, 40\penalty0 (1):\penalty0 120--145,
  2011.
\newblock \doi{10.1007/s10851-010-0251-1}.
\newblock URL \url{https://doi.org/10.1007/s10851-010-0251-1}.

\bibitem[Applegate et~al.(2021{\natexlab{b}})Applegate, D{\'\i}az, Hinder, Lu,
  Lubin, O'Donoghue, and Schudy]{PDLP-2021}
David Applegate, Mateo D{\'\i}az, Oliver Hinder, Haihao Lu, Miles Lubin,
  Brendan O'Donoghue, and Warren Schudy.
\newblock Practical large-scale linear programming using primal-dual hybrid
  gradient.
\newblock \emph{Advances in Neural Information Processing Systems},
  34:\penalty0 20243--20257, 2021{\natexlab{b}}.

\bibitem[Stellato et~al.(2020)Stellato, Banjac, Goulart, Bemporad, and
  Boyd]{stellato2020osqp}
Bartolomeo Stellato, Goran Banjac, Paul Goulart, Alberto Bemporad, and Stephen
  Boyd.
\newblock Osqp: An operator splitting solver for quadratic programs.
\newblock \emph{Mathematical Programming Computation}, 12\penalty0
  (4):\penalty0 637--672, 2020.

\bibitem[Garstka et~al.(2021{\natexlab{a}})Garstka, Cannon, and
  Goulart]{garstka2021cosmo}
Michael Garstka, Mark Cannon, and Paul Goulart.
\newblock Cosmo: A conic operator splitting method for convex conic problems.
\newblock \emph{Journal of Optimization Theory and Applications}, 190\penalty0
  (3):\penalty0 779--810, 2021{\natexlab{a}}.

\bibitem[O'Donoghue(2021)]{o2021operator}
Brendan O'Donoghue.
\newblock Operator splitting for a homogeneous embedding of the linear
  complementarity problem.
\newblock \emph{SIAM Journal on Optimization}, 31\penalty0 (3):\penalty0
  1999--2023, 2021.

\bibitem[O'Donoghue et~al.(2023)O'Donoghue, Chu, Parikh, and Boyd]{scs}
Brendan O'Donoghue, Eric Chu, Neal Parikh, and Stephen Boyd.
\newblock {SCS}: Splitting conic solver, version 3.2.4.
\newblock \url{https://github.com/cvxgrp/scs}, November 2023.

\bibitem[Cortes and Vapnik(1995)]{cortes1995support}
Corinna Cortes and Vladimir Vapnik.
\newblock Support-vector networks.
\newblock \emph{Machine learning}, 20\penalty0 (3):\penalty0 273--297, 1995.

\bibitem[Tibshirani(1996)]{tibshirani1996regression}
Robert Tibshirani.
\newblock Regression shrinkage and selection via the lasso.
\newblock \emph{Journal of the Royal Statistical Society: Series B
  (Methodological)}, 58\penalty0 (1):\penalty0 267--288, 1996.

\bibitem[Markowitz(1991)]{markowitz1991foundations}
Harry~M Markowitz.
\newblock Foundations of portfolio theory.
\newblock \emph{The journal of finance}, 46\penalty0 (2):\penalty0 469--477,
  1991.

\bibitem[Luo et~al.(1997)Luo, Sturm, and Zhang]{Luo-1997}
Z-Q. Luo, J.F. Sturm, and Shuzhong Zhang.
\newblock Duality results for conic convex programming.
\newblock Econometric Institute Research Papers EI 9719/A, Erasmus University
  Rotterdam, Erasmus School of Economics (ESE), Econometric Institute, 1997.

\bibitem[{Zhang}(2004)]{Zhang2004}
Shuzhong {Zhang}.
\newblock A new self-dual embedding method for convex programming.
\newblock \emph{Journal of Global Optimization}, 29\penalty0 (4):\penalty0
  479--496, 2004.

\bibitem[Banjac et~al.(2019)Banjac, Goulart, Stellato, and
  Boyd]{banjac2019infeasibility}
Goran Banjac, Paul Goulart, Bartolomeo Stellato, and Stephen Boyd.
\newblock Infeasibility detection in the alternating direction method of
  multipliers for convex optimization.
\newblock \emph{Journal of Optimization Theory and Applications}, 183:\penalty0
  490--519, 2019.

\bibitem[Applegate et~al.(2021{\natexlab{c}})Applegate, D{\'\i}az, Lu, and
  Lubin]{applegate2021infeasibility}
David Applegate, Mateo D{\'\i}az, Haihao Lu, and Miles Lubin.
\newblock Infeasibility detection with primal-dual hybrid gradient for
  large-scale linear programming.
\newblock \emph{arXiv preprint arXiv:2102.04592}, 2021{\natexlab{c}}.

\bibitem[Wright et~al.(1999)Wright, Nocedal, et~al.]{wright1999numerical}
Stephen Wright, Jorge Nocedal, et~al.
\newblock Numerical optimization.
\newblock \emph{Springer Science}, 35\penalty0 (67-68):\penalty0 7, 1999.

\bibitem[W{\"a}chter and Biegler(2006)]{wachter2006implementation}
Andreas W{\"a}chter and Lorenz~T Biegler.
\newblock On the implementation of an interior-point filter line-search
  algorithm for large-scale nonlinear programming.
\newblock \emph{Mathematical programming}, 106\penalty0 (1):\penalty0 25--57,
  2006.

\bibitem[Vanderbei(1999)]{vanderbei1999loqo}
Robert~J Vanderbei.
\newblock Loqo: An interior point code for quadratic programming.
\newblock \emph{Optimization methods and software}, 11\penalty0 (1-4):\penalty0
  451--484, 1999.

\bibitem[Pokutta(2020)]{Pokutta-2020}
Sebastian Pokutta.
\newblock Restarting algorithms: sometimes there is free lunch.
\newblock In \emph{International Conference on Integration of Constraint
  Programming, Artificial Intelligence, and Operations Research}, pages 22--38.
  Springer, 2020.

\bibitem[Gleixner et~al.(2022)Gleixner, Gottwald, and Hoen]{Papilo}
Ambros Gleixner, Leona Gottwald, and Alexander Hoen.
\newblock Papilo: A parallel presolving library for integer and linear
  programming with multiprecision support.
\newblock \emph{arXiv preprint arXiv:2206.10709}, 2022.

\bibitem[Pock and Chambolle(2011)]{PockChambolle-2011}
Thomas Pock and Antonin Chambolle.
\newblock Diagonal preconditioning for first order primal-dual algorithms in
  convex optimization.
\newblock In \emph{2011 International Conference on Computer Vision}, pages
  1762--1769. IEEE, 2011.

\bibitem[Ruiz(2001)]{Ruiz-2001}
Daniel Ruiz.
\newblock A scaling algorithm to equilibrate both rows and columns norms in
  matrices.
\newblock Technical report, CM-P00040415, 2001.

\bibitem[Nesterov(2014)]{nesterov-subgradient}
Yu~Nesterov.
\newblock Subgradient methods for huge-scale optimization problems.
\newblock \emph{Mathematical Programming}, 146\penalty0 (1):\penalty0 275--297,
  2014.

\bibitem[Gleixner et~al.(2021)Gleixner, Hendel, Gamrath, Achterberg, Bastubbe,
  Berthold, Christophel, Jarck, Koch, Linderoth, et~al.]{gleixner2021miplib}
Ambros Gleixner, Gregor Hendel, Gerald Gamrath, Tobias Achterberg, Michael
  Bastubbe, Timo Berthold, Philipp Christophel, Kati Jarck, Thorsten Koch, Jeff
  Linderoth, et~al.
\newblock Miplib 2017: data-driven compilation of the 6th mixed-integer
  programming library.
\newblock \emph{Mathematical Programming Computation}, 13\penalty0
  (3):\penalty0 443--490, 2021.

\bibitem[Koch(2003)]{koch2003final}
Thorsten Koch.
\newblock The final netlib-lp results.
\newblock 2003.

\bibitem[Davis et~al.(2013)Davis, Duff, Amestoy, Gilbert, Larimore, Natarajan,
  Chen, Hager, and Rajamanickam]{suitesparse}
Tim Davis, I~Duff, P~Amestoy, J~Gilbert, S~Larimore, E~Palamadai Natarajan,
  Y~Chen, W~Hager, and S~Rajamanickam.
\newblock Suitesparse: A suite of sparse matrix packages.
\newblock \url{https://github.com/DrTimothyAldenDavis/SuiteSparse}, 2013.
\newblock {Accessed}: March 15, 2024.

\bibitem[Chang and Lin(2011)]{chang2011libsvm}
Chih-Chung Chang and Chih-Jen Lin.
\newblock Libsvm: a library for support vector machines.
\newblock \emph{ACM transactions on intelligent systems and technology (TIST)},
  2\penalty0 (3):\penalty0 1--27, 2011.

\bibitem[Andersen and Ye(1999)]{andersen1999homogeneous}
Erling~D Andersen and Yinyu Ye.
\newblock On a homogeneous algorithm for the monotone complementarity problem.
\newblock \emph{Mathematical Programming}, 84\penalty0 (2):\penalty0 375--399,
  1999.

\bibitem[Friberg(2016)]{friberg2016cblib}
Henrik~A Friberg.
\newblock Cblib 2014: a benchmark library for conic mixed-integer and
  continuous optimization.
\newblock \emph{Mathematical Programming Computation}, 8\penalty0 (2):\penalty0
  191--214, 2016.

\bibitem[Garstka et~al.(2021{\natexlab{b}})Garstka, Cannon, and
  Goulart]{Garstka_2021}
Michael Garstka, Mark Cannon, and Paul Goulart.
\newblock {COSMO}: A conic operator splitting method for convex conic problems.
\newblock \emph{Journal of Optimization Theory and Applications}, 190\penalty0
  (3):\penalty0 779--810, 2021{\natexlab{b}}.
\newblock \doi{10.1007/s10957-021-01896-x}.
\newblock URL \url{https://doi.org/10.1007/s10957-021-01896-x}.

\bibitem[Combettes and Wajs(2005)]{combettes2005signal}
Patrick~L Combettes and Val{\'e}rie~R Wajs.
\newblock Signal recovery by proximal forward-backward splitting.
\newblock \emph{Multiscale modeling \& simulation}, 4\penalty0 (4):\penalty0
  1168--1200, 2005.

\end{thebibliography}
\clearpage
\appendix

\section{Proof of Theorem~\ref{thm-pq}}

\label{proof-thm-pq}

	We shall prove the result by induction. (i) At iteration $j = 0$,  the result holds true from our definition. (ii) Assuming it holds true for iteration $j = i$, we will prove that the result still holds true for iteration $j = i+1$ in two steps:
	
	\textit{Step 1}: We claim that
	\begin{equation}
		\bu^{(k)}_{i} + \bv^{(k)}_{i} = \tilde{\bu}^{(k)}_{i+1} + \tilde{\bv}^{(k)}_{i+1} \label{step1}.
	\end{equation}
	Indeed, \eqref{admm1} can be rewritten as
	\begin{equation}
		(\tilde{\bu}^{(k)}_{i+1}, \tilde{\bv}^{(k)}_{i+1}) = \Pi_{\mathcal{P}} (\bu^{(k)}_{i}+\bv^{(k)}_{i},\bu^{(k)}_{i}+\bv^{(k)}_{i}), \label{p1}
	\end{equation}
	where $\mathcal{P} = \{(\bu,\bv): \bQ\bu=\bv\}$. Moreover, since $\bQ$ is skew-symmetric, the orthogonal complement of $\mathcal{P}$ is $\mathcal{P}^{\perp} = \{(\bv,\bu):\bQ\bu=\bv\}$. Therefore, we conclude that
	\begin{equation*}
		(\bu,\bv) = \Pi_{\mathcal{P}}(\bz,\bz) \quad \text{if and only if}\quad (\bv,\bu) = \Pi_{\mathcal{P}^{\perp}}(\bz,\bz),
	\end{equation*}
	because the two projections are identical for reversed output arguments. This implies that
	\begin{equation}
		(\tilde{\bv}^{(k)}_{i+1}, \tilde{\bu}^{(k)}_{i+1}) = \Pi_{\mathcal{P}^{\perp}} (\bu^{(k)}_{i}+\bv^{(k)}_{i},\bu^{(k)}_{i}+\bv^{(k)}_{i}). \label{p2}
	\end{equation}
	Then we combine (\ref{p1}) and (\ref{p2}) to get the desired result.
	
	\textit{Step 2}: 	Given $	\bp^{(k)}_{i}=\bv^{(k)}_{i}$, $\bq^{(k)}_{i}=\bu^{(k)}_{i}$ and $\xi^{(k)}_{i} = -\bx_{0}^{T}\bs_{0}-1$, our goal is to show that 
	\begin{equation*}
		\bp^{(k)}_{i+1}=\bv^{(k)}_{i+1},\quad \bq^{(k)}_{i+1}=\bu^{(k)}_{i+1}.
	\end{equation*}
Note that \eqref{step1} can be rewritten as 
	\begin{align*}
		 & \tilde{\bu}^{(k)}_{i+1} - \bv^{(k)}_{i} = - (\tilde{\bv}^{(k)}_{i+1} - \bu^{(k)}_{i}).
	\end{align*}
	We denote $\sigma = \tilde{\bu}^{(k)}_{i+1} - \bv^{(k)}_{i}$. It is easy to verify that functions in \eqref{eq:B-conic} form pairs of the function $h(t)$ and $w(t)$ satisfying $w(t) = h^{*}(-t)$, where $h^*$ denotes the conjugate function of $h$. From Moreau's decomposition (see Section 2.5 in \cite{combettes2005signal} for more details), for such pairs of functions $(h(t),w(t))$ in $B$, we have
	\begin{equation*}
	    \bu^{(k)}_{i+1} - \bv^{(k)}_{i+1} = \text{prox}_{h}(\sigma) - \text{prox}_{w}(-\sigma) = \sigma.
	\end{equation*}
	Then from \eqref{admm3} we have
	\begin{align*}
		\bp^{(k)}_{i+1} &= \bv^{(k)}_{i} + \bu^{(k)}_{i+1} - \tilde{\bu}^{(k)}_{i+1} = \bu^{(k)}_{i+1} - \sigma = \bv^{(k)}_{i+1}, \\
		\bq^{(k)}_{i+1} &= \bu^{(k)}_{i} + \bv^{(k)}_{i+1} - \tilde{\bv}^{(k)}_{i+1} = \bv^{(k)}_{i+1} + \sigma = \bu^{(k)}_{i+1}.
	\end{align*}

	This completes the proof.
 
 \hfill$\square$

\section{Proof of Proposition~\ref{prop:complexity-admm}}
\label{proof-prop-5}
	Observing that
	\begin{align}
		&\nabla^2_{x} B(\bu,\bv,\mu^{(k)})\succeq \dfrac{\mu^{(k)}}{C_D}\bI\label{convexity1}\\
		&\nabla^2_{s} B(\bu,\bv,\mu^{(k)})\succeq \dfrac{\mu^{(k)}}{C_D}\bI\label{convexity2}\\
		&\nabla^2_\tau B(\bu,\bv,\mu^{(k)})\geq \dfrac{\mu^{(k)}}{D}\label{convexity3}\\
		&\nabla^2_\kappa B(\bu,\bv,\mu^{(k)})\geq \dfrac{\mu^{(k)}}{D}\label{convexity4}
	\end{align}
	For convenience, we assume that $C_D\geq D$. The optimality condition of problem (\ref{conic-subproblem}) implies that
	\begin{equation*}
		\beta(\bp_k^*,\bq_k^*)\in\partial \mathds{1}(\bQ\bu=\bv)[\bu_k^*,\bv_k^*],\quad-\beta(\bp_k^*,\bq_k^*)\in\partial B(\bu_k^*,\bv_k^*,\mu^{(k)}).
	\end{equation*}
	Using the convexity of the two functions and (\ref{convexity1})(\ref{convexity2})(\ref{convexity3})(\ref{convexity4}) we have
	\begin{align}
		&0\leq\beta(\bu_k^*-\tilde{\bu}_{i+1}^{(k)},\bv_k^*-\tilde{\bv}_{i+1}^{(k)})^\top (\bp_k^*,\bq_k^*)=\beta(\bu_k^*-\tilde{\bu}_{i+1}^{(k)},\bv_k^*-\tilde{\bv}_{i+1}^{(k)})^\top(\bv_k^*,\bu_k^*)\label{58}\\
	    \begin{split}
			&B(\bu_k^*,\bv_k^*,\mu^{(k)})-B(\bu_{i+1}^{(k)}-\bv_{i+1}^{(k)},\mu^{(k)})+\dfrac{\mu^{(k)}}{2C_D}(||\bx_{i+1}^{(k)}-\bx_k^*||^2+||\bs_{i+1}^{(k)}-\bs_k^*||^2+(\tau_{i+1}^{(k)}-\tau_k^*)^2+(\kappa_{i+1}^{(k)}-\kappa_k^*)^2)\\
			\leq&-\beta(\bu_k^*-\bu_{i+1}^{(k)},\bv_k^*-\bv_{i+1}^{(k)})^\top (\bp_k^*,\bq_k^*)
			=-\beta(\bu_k^*-\bu_{i+1}^{(k)},\bv_k^*-\bv_{i+1}^{(k)})^\top (\bv_k^*,\bu_k^*).
		\end{split}\label{59}
	\end{align}
	Summing up (\ref{58})(\ref{59}) we have
	\begin{equation}\label{60}
		\begin{split}
			&B(\bu_k^*,\bv_k^*,\mu^{(k)})-B(\bu_{i+1}^{(k)}-\bv_{i+1}^{(k)},\mu^{(k)})+\dfrac{\mu^{(k)}}{2C_D}(||\bx_{i+1}^{(k)}-\bx_k^*||^2+||\bs_{i+1}^{(k)}-\bs_k^*||^2+(\tau_{i+1}^{(k)}-\tau_k^*)^2+(\kappa_{i+1}^{(k)}-\kappa_k^*)^2)\\
			\leq {}&\beta(\bu_{i+1}^{(k)}-\tilde{\bu}_{i+1}^{(k)},\bv_{i+1}^{(k)}-\tilde{\bv}_{i+1}^{(k)})^\top (\bv_k^*,\bu_k^*)
		\end{split}
	\end{equation}
	With the definition of $\bu_{i+1}^{(k)},\bv_{i+1}^{(k)},\tilde{\bu}_{i+1}^{(k)},\tilde{\bv}_{i+1}^{(k)}$, we have
	\begin{align*}
		&\beta(\bp_{i+1}^{(k)},\bq_{i+1}^{(k)})-\beta(\bu_{i+1}^{(k)},\bv_{i+1}^{(k)})+\beta(\bu_i^{(k)},\bv_i^{(k)})\leq \partial \mathds{1}(\bQ\bu=\bv)[\tilde{\bu}_{i+1}^{(k)},\tilde{\bv}_{i+1}^{(k)}]\\
		&-\beta(\bp_{i+1}^{(k)},\bq_{i+1}^{(k)})\in \partial B(\bu_{i+1}^{(k)},\bv_{i+1}^{(k)},\mu^{(k)}).
	\end{align*}
	Therefore
	\begin{align}
		\begin{split}\label{61}
			0\leq {}&\beta(\tilde{\bu}_{i+1}^{(k)}-\bu_k^*,\tilde{\bv}_{i+1}^{(k)},\bv_k^*)^\top (\bp_{i+1}^{(k)}-\bu_{i+1}^{(k)}+\bu_i^{(k)},\bq_{i+1}^{(k)}-\bv_{i+1}^{(k)}+\bv_i^{(k)})\\
			={}&\beta(\tilde{\bu}_{i+1}^{(k)}-\bu_k^*,\tilde{\bv}_{i+1}^{(k)},\bv_k^*)^\top (\bp_{i+1}^{(k)},\bq_{i+1}^{(k)})+\dfrac{\beta}{2}(||\bu_i^{(k)}-\bu_k^*||^2+||\bv_i^{(k)}-\bv_k^*||^2-||\bu_i^{(k)}-\tilde{\bu}_{i+1}^{(k)}||^2-||\bv_i^{(k)}-\tilde{\bv}_{i+1}^{(k)}||^2)\\
			&-\dfrac{\beta}{2}(||\bu_{i+1}^{(k)}-\bu_k^*||^2+||\bv_{i+1}^{(k)}-\bv_k^*||^2-||\bu_{i+1}^{(k)}-\tilde{\bu}_{i+1}^{(k)}||^2-||\bv_{i+1}^{(k)}-\tilde{\bv}_{i+1}^{(k)}||^2)
		\end{split}\\
		\begin{split}\label{62}
			&B(\bu_{i+1}^{(k)},\bv_{i+1}^{(k)},\mu^{(k)})-B(\bu_k^*,\bv_k^*,\mu^{(k)})+\dfrac{\mu^{(k)}}{2C_D}(||\bx_{i+1}^{(k)}-\bx_k^*||^2+||\bs_{i+1}^{(k)}-\bs_k^*||^2+(\tau_{i+1}^{(k)}-\tau_k^*)^2+(\kappa_{i+1}^{(k)}-\kappa_k^*)^2)\\
			\leq &-\beta(\bu_{i+1}^{(k)}-\bu_k^*,\bv_{i+1}^{(k)}-\bv_k^*)^\top (\bp_{i+1}^{(k)},\bq_{i+1}^{(k)}).
		\end{split}
	\end{align}
	Summing up (\ref{61}) and (\ref{62}) gives
	\begin{equation}\label{63}
		\begin{split}
			&B(\bu_{i+1}^{(k)},\bv_{i+1}^{(k)},\mu^{(k)})-B(\bu_k^*,\bv_k^*,\mu^{(k)})+\dfrac{\mu^{(k)}}{4D}(||\bx_{i+1}^{(k)}-\bx_k^*||^2+||\bs_{i+1}^{(k)}-\bs_k^*||^2+(\tau_{i+1}^{(k)}-\tau_k^*)^2+(\kappa_{i+1}^{(k)}-\kappa_k^*)^2)\\
			\leq {} &\beta(\tilde{\bu}_{i+1}^{(k)}-\bu_{i+1}^{(k)},\tilde{\bv}_{i+1}^{(k)}-\bv_{i+1}^{(k)})^\top(\bv_{i+1}^{(k)},\bu_{i+1}^{(k)})+\dfrac{\beta}{2}(||\bu_i^{(k)}-\bu_k^*||^2+||\bv_i^{(k)}-\bv_k^*||^2-||\bu_i^{(k)}-\tilde{\bu}_{i+1}^{(k)}||^2-||\bv_i^{(k)}-\tilde{\bv}_{i+1}^{(k)}||^2)\\
			&-\dfrac{\beta}{2}(||\bu_{i+1}^{(k)}-\bu_k^*||^2+||\bv_{i+1}^{(k)}-\bv_k^*||^2-||\bu_{i+1}^{(k)}-\tilde{\bu}_{i+1}^{(k)}||^2-||\bv_{i+1}^{(k)}-\tilde{\bv}_{i+1}^{(k)}||^2).
		\end{split}
	\end{equation}
	Observing that
	\begin{equation}\label{64}
		\begin{split}
			&(\bu_{i+1}^{(k)}-\tilde{\bu}_{i+1}^{(k)})^\top (\bv_k^*,\bu_k^*)+(\tilde{\bu}_{i+1}^{(k)}-\bu_{i+1}^{(k)},\tilde{\bv}_{i+1}^{(k)}-\bv_{i+1}^{(k)})^\top (\bv_{i+1}^{(k)},\bu_{i+1}^{(k)})\\
			={}&(\tilde{\bu}_{i+1}^{(k)}-\bu_{i+1}^{(k)},\tilde{\bv}_{i+1}^{(k)}-\bv_{i+1}^{(k)})^\top (\bv_{i+1}^{(k)}-\bv_k^*,\bu_{i+1}^{(k)}-\bu_k^*)\\
			={}&(\bp_i^{(k)}-\bp_{i+1}^{(k)},\bq_i^{(k)}-\bq_{i+1}^{(k)})^\top (\bv_{i+1}^{(k)}-\bv_k^*,\bu_{i+1}^{(k)}-\bu_k^*)\\
			={}&(\bv_i^{(k)}-\bv_{i+1}^{(k)},\bu_i^{(k)}-\bu_{i+1}^{(k)})^\top (\bv_{i+1}^{(k)}-\bv_k^*,\bu_{i+1}^{(k)}-\bu_k^*)\\
			= {}&\dfrac{1}{2}(||\bu_i^{(k)}-\bu_k^*||^2+||\bv_i^{(k)}-\bv_k^*||^2-||\bu_{i+1}^{(k)}-\bu_k^*||^2-||\bv_{i+1}^{(k)}-\bv_k^*||^2-||\bu_i^{(k)}-\bu_{i+1}^{(k)}||^2-||\bv_i^{(k)}-\bv_{i+1}^{(k)}||^2).
		\end{split}
	\end{equation}
	Summing (\ref{60}), (\ref{63}), and using (\ref{64}), we have
	\begin{equation}\label{67}
		\begin{split}
			&\dfrac{1}{2}(||\bu_i^{(k)}-\tilde{\bu}_{i+1}||^2+||\bv_i^{(k)}-\tilde{\bv}_{i+1}^{(k)}||^2)+\dfrac{\mu^{(k)}}{2C_D\beta}(||\bx_{i+1}^{(k)}-\bx_k^*||^2+||\bs_{i+1}^{(k)}-\bs_k^*||^2+(\tau_{i+1}^{(k)}-\tau_k^*)^2+(\kappa_{i+1}^{(k)}-\kappa_k^*)^2)\\
			\leq {} &||\bu_i^{(k)}-\bu_k^*||^2+||\bv_i^{(k)}-\bv_k^*||^2-||\bu_{i+1}^{(k)}-\bu_k^*||^2-||\bv_{i+1}^{(k)}-\bv_k^*||^2.
		\end{split}
	\end{equation}
	Using the optimality condition of (\ref{admm2}) we have
	\begin{equation}\label{65}
		0\leq B(\bu,\bv,\mu^{(k)})-B(\bu_{i+1}^{(k)},\bv_{i+1}^{(k)},\mu^{(k)})+\beta(\bu-\bu_{i+1}^{(k)},\bv-\bv_{i+1}^{(k)})^\top (\bp_{i+1}^{(k)},\bq_{i+1}^{(k)})
	\end{equation}
	and
	\begin{equation}\label{66}
		0\leq B(\bu,\bv,\mu^{(k)})-B(\bu_{i+1}^{(k)},\bv_{i+1}^{(k)},\mu^{(k)})+\beta(\bu-\bu_{i+1}^{(k)},\bv-\bv_{i+1}^{(k)})^\top (\bp_{i+1}^{(k)},\bq_{i+1}^{(k)}).
	\end{equation}
	Let $(\bu,\bv)=(\bu_i^{(k)},\bv_i^{(k)})$ in (\ref{65}) and $(\bu,\bv)=(\bu_{i+1},\bv_{i+1}^{(k)})$ in (\ref{66}), and sum up the two inequalities we have
	\begin{equation}
		\begin{aligned}
			0 &\le-(\bu_i^{(k)}-\bu_{i+1}^{(k)},\bv_i^{(k)}-\bv_{i+1}^{(k)})^\top (\bp_i^{(k)}-\bp_{i+1}^{(k)},\bq_i^{(k)}-\bq_{i+1}^{(k)})\\
			&=-(\bu_i^{(k)}-\bu_{i+1}^{(k)},\bv_i^{(k)}-\bv_{i+1}^{(k)})^\top (\bu_{i+1}^{(k)}-\tilde{\bu}_{i+1}^{(k)},\bv_i^{(k)}-\tilde{\bv}_{i+1}^{(k)})\\
			&=\dfrac{1}{2}\left(||\bu_i^{(k)}-\tilde{\bu}_{i+1}^{(k)}||^2-||\bu_{i+1}^{(k)}-\tilde{\bu}_{i+1}^{(k)}||^2-||\bu_i^{(k)}-\bu_{i+1}^{(k)}||^2\right)\\
			&\quad +\dfrac{1}{2}\left(||\bv_i^{(k)}-\tilde{\bv}_{i+1}^{(k)}||^2-||\bv_{i+1}^{(k)}-\tilde{\bv}_{i+1}^{(k)}||^2-||\bv_i^{(k)}-\bv_{i+1}^{(k)}||^2\right),
		\end{aligned}
	\end{equation}
	which implies that
	\begin{equation}\label{68}
		||\bu_{i+1}^{(k)}-\tilde{\bu}_{i+1}^{(k)}||^2+||\bv_{i+1}^{(k)}-\tilde{\bv}_{i+1}^{(k)}||^2\leq ||\bu_{i}^{(k)}-\tilde{\bu}_{i+1}^{(k)}||^2+||\bv_{i}^{(k)}-\tilde{\bv}_{i+1}^{(k)}||^2.
	\end{equation}
	Summing up (\ref{67}) and (\ref{68}) we have
	\begin{equation}
		\begin{split}
			&\dfrac{1}{2}\left(||\bu_{i+1}^{(k)}-\tilde{\bu}_{i+1}^{(k)}||^2+||\bv_{i+1}^{(k)}-\tilde{\bv}_{i+1}^{(k)}||^2\right)\\
			&+\dfrac{\mu^{(k)}}{2C_D\beta}(||\bx_{i+1}^{(k)}-\bx_k^*||^2+||\bs_{i+1}^{(k)}-\bs_k^*||^2+(\tau_{i+1}^{(k)}-\tau_k^*)^2+(\kappa_{i+1}^{(k)}-\kappa_k^*)^2)\\
			\leq{}& ||\bu_i^{(k)}-\bu_k^*||^2+||\bv_i^{(k)}-\bv_k^*||^2-||\bu_{i+1}^{(k)}-\bu_k^*||^2-||\bv_{i+1}^{(k)}-\bv_k^*||^2.
		\end{split}
	\end{equation}
	Using the linear constraint of (\ref{conic-subproblem}) we have (by denoting $C_4=6\lambda_\text{max}(\bA^\top \bA)\max\{1,||c||^2,||r_d||^2\}/\lambda^2_\text{min}(\bA\bA^\top)$)
	\begin{equation}
		\begin{split}\label{69}
			&||\by_{i+1}^{(k)}-\by_k^*||^2\\
			\leq{} &2\left(||\by_{i+1}^{(k)}-\tilde{\by}_{i+1}^{(k)}||^2+||\tilde{\by}_{i+1}^{(k)}-\by_k^*||^2\right)\\
			={}&2||\by_{i+1}^{(k)}-\tilde{\by}_{i+1}^{(k)}||^2+2||(\bA\bA^\top)^{-1}\bA(\bA^\top\tilde{\by}_{i+1}^{(k)}-\bA^T y_k^*)||^2\\
			\leq{}&2||\by_{i+1}^{(k)}-\tilde{\by}_{i+1}^{(k)}||^2+C_4\left(||\tilde{\by}_{i+1}^{(k)}-\bs_k^*||^2+(\tilde{\by}_{i+1}^{(k)}-\tau_k^*)^2+(\tilde{\theta}_{i+1}^{(k)}-\theta_k^*)^2\right)\\
			\leq{}&2C_4\left(||\tilde{\by}_{i+1}^{(k)}-\bs_{i+1}^{(k)}||^2+||\bs_{i+1}^{(k)}-\bs_k^*||^2+(\tilde{\by}_{i+1}^{(k)}-\tau_{i+1}^{(k)})^2+(\tau_{i+1}^{(k)}-\tau_k^*)^2\right)\\
			{}&+C_4(\theta_{i+1}^{(k)}-\theta_k^*)^2+2||\by_{i+1}^{(k)}-\tilde{\by}_{i+1}^{(k)}||^2
		\end{split}
	\end{equation}
	and (by denoting $C_5=3\max\{||\bb||^2,||\bA||^2/||\br_p||^2\}$)
	\begin{equation}
		\begin{split}\label{70}
			{}&(\theta_{i+1}^{(k)}-\theta_k^*)^2\\
			={}&(\tilde{\theta}_{i+1}^{(k)}-\theta_k^*)^2=\dfrac{1}{||\br_p||^2}||\bb\tilde{\by}_{i+1}^{(k)}-\bA\tilde{\bx}_{i+1}^{(k)}+\tilde{\br}_{i+1}^{(k)}-\bb\tau_k^*+\bA\bx_k^*-\br_k^*||^2\\
			\leq{} &C_5(||\tilde{\bx}_{i+1}^{(k)}-\bx_k^*||^2+(\tilde{\by}_{i+1}^{(k)}-\tau_k^*)^2+||\tilde{\br}_{i+1}^{(k)}-\br_k^*||^2)\\
			\leq{} &2C_5(||\bx_{i+1}^{(k)}-\tilde{\bx}_{i+1}^{(k)}||^2+||\bx_{i+1}^{(k)}-\bx_k^*||^2+(\tau_{i+1}^{(k)}-\tilde{\by}_{i+1}^{(k)})^2+(\tau_{i+1}^{(k)}-\tau_k^*)^2)+C_5||\br_{i+1}^{(k)}-\tilde{\br}_k^*||^2.
		\end{split}
	\end{equation}
	Summing up (\ref{69}) and (\ref{70}) we have
	\begin{equation}
		\begin{split}
			&||\by_{i+1}^{(k)}-\by_k^*||^2+(\theta_{i+1}^{(k)}-\theta_k^*)^2\\
			\leq {} &2(C_4+C_5)\left(||\bu_{i+1}^{(k)}-\tilde{\bu}_{i+1}^{(k)}||^2+||\bv_{i+1}^{(k)}-\tilde{\bv}_{i+1}^{(k)}||^2\right)\\
			{}&+2C_4\left(||\tilde{\by}_{i+1}^{(k)}-\bs_{i+1}^{(k)}||^2+(\tau_{i+1}^{(k)}-\tau_k^*)^2+(\theta_{i+1}^{(k)}-\theta_k^*)^2\right) \\
			{} &+2C_5\left(||\bx_{i+1}^{(k)}-\bx_k^*||^2+(\tau_{i+1}^{(k)}-\tau_k^*)^2\right)\\
			\leq {} &2(C_4+C_5)\left(||\bu_{i+1}^{(k)}-\tilde{\bu}_{i+1}^{(k)}||^2+||\bv_{i+1}^{(k)}-\tilde{\bv}_{i+1}^{(k)}||^2\right)\\
			{} &+2C_4\left(||\tilde{\by}_{i+1}^{(k)}-\bs_{i+1}^{(k)}||^2+(\tau_{i+1}^{(k)}-\tau_k^*)^2\right)+2C_5\left(||\bx_{i+1}^{(k)}-\bx_k^*||^2+(\tau_{i+1}^{(k)}-\tau_k^*)^2\right)\\
			{} &+4C_4C_5(||\bx_{i+1}^{(k)}-\tilde{\bx}_{i+1}^{(k)}||^2+||\bx_{i+1}^{(k)}-\bx_k^*||^2+(\tau_{i+1}^{(k)}-\tilde{\by}_{i+1}^{(k)})^2+(\tau_{i+1}^{(k)}-\tau_k^*)^2)\\
			{}&+2C_4C_5||\tilde{\br}_{i+1}^{(k)}-\br_{i+1}^{(k)}||^2\\
			\leq {} &C_3(||\bu_{i+1}^{(k)}-\tilde{\bu}_{i+1}^{(k)}||^2)+C_3(||\bx_{i+1}^{(k)}-\bx_k^*||^2+||\bs_{i+1}^{(k)}-\bs_k^*||^2+(\tau_{i+1}^{(k)}-\tau_k^*)^2).
		\end{split}
	\end{equation}
	Finally, we have
	\begin{equation}
		\begin{split}
			&\min\left\{\dfrac{1}{C_3},\dfrac{\mu^{(k)}}{2C_DC_3\beta}\right\}(||\bu_{i+1}^{(k)}-\bu_k^*||^2+||\bv_{i+1}^{(k)}-\bv_k^*||^2)\\
			\leq{} &\dfrac{\mu^{(k)}}{2C_DC_3\beta}(||\bx_{i+1}^{(k)}-\bx_k^*||^2+||\bs_{i+1}^{(k)}-\bs_k^*||^2+(\tau_{i+1}^{(k)}-\tau_k^*)^2+(\kappa_{i+1}^{(k)}-\kappa_k^*)^2)\\
			&+\min\{\dfrac{1}{2C_3},\dfrac{\mu^{(k)}}{2C_DC_3\beta}\}(||\by_{i+1}^{(k)}-\by_k^*||^2+(\theta_{i+1}^{(k)}-\theta_k^*)^2)\\
			\leq {} &\dfrac{\mu^{(k)}}{2C_D\beta}(||\bx_{i+1}^{(k)}-\bx_k^*||^2+||\bs_{i+1}^{(k)}-\bs_k^*||^2+(\tau_{i+1}^{(k)}-\tau_k^*)^2+(\kappa_{i+1}^{(k)}-\kappa_k^*)^2)\\
			&+\dfrac{1}{2}(||\bu_{i+1}^{(k)}-\tilde{\bu}_{i+1}^{(k)}||^2+||\bv_{i+1}^{(k)}-\tilde{\bv}_{i+1}^{(k)}||^2)\\
			&+\dfrac{\mu^{(k)}}{2C_D\beta}(||\bx_{i+1}^{(k)}-\bx_k^*||^2+||\bs_{i+1}^{(k)}-\bs_k^*||^2+(\tau_{i+1}^{(k)}-\tau_k^*)^2)\\
			\leq {}&+\dfrac{\mu^{(k)}}{2C_D\beta}\big[||\bx_{i+1}^{(k)}-\bx_k^*||^2+||\bs_{i+1}^{(k)}-s_k^*||^2+(\tau_{i+1}^{(k)}-\tau_k^*)^2+(\kappa_{i+1}^{(k)}-\kappa_k^*)^2\big]\\
			{} & + \dfrac{1}{2}(||\bu_i^{(k)}-\tilde{\bu}_{i+1}||^2+||\bv_i^{(k)}-\tilde{\bv}_{i+1}^{(k)}||^2)\\
			\leq {}&||\bu_i^{(k)}-\bu_k^*||^2+||\bv_i^{(k)}-\bv_k^*||^2-||\bu_{i+1}^{(k)}-\bu_k^*||^2-||\bv_{i+1}^{(k)}-\bv_k^*||^2.
		\end{split}
	\end{equation}
	Therefore, we have
	\begin{equation}
		\begin{split}
			||\bu_{N_k}^{(k)}-\bu_k^*||^2+||\bv_{N_k}^{(k)}-\bv_k^*||^2&\leq \left(1+\min\left\{\dfrac{1}{C_3},\dfrac{\mu^{(k)}}{2C_DC_3\beta}\right\}\right)^{-N_k}(||\bu_0^{(k)}-\bu_k^*||^2+||\bv_0^{(k)}-\bv_k^*||^2)\\
			&\leq 2C_1\left(1+\min\left\{\dfrac{1}{C_3},\dfrac{\mu^{(k)}}{2C_DC_3\beta}\right\}\right)^{-N_k}.
		\end{split}
	\end{equation}
	Combining the termination criteria and the following inequality
	\begin{equation}
		||\bQ \bu_{N_k}^{(k)}-\bv_{N_k}^{(k)}||^2\leq (2+2||\bQ||^2)(||\bu_{N_k}^{(k)}-\bu_k^*||^2+||\bv_{N_k}^{(k)}-\bv_k^*||^2),
	\end{equation}
	we have
	\begin{equation}
		4C_1(1+||\bQ||^2)\left(1+\min\left\{\dfrac{1}{C_3},\dfrac{\mu^{(k)}}{2C_DC_3\beta}\right\}\right)^{-N_k}\geq \mu^{(k)},
	\end{equation}
	which completes our proof.

  \hfill$\square$

\section{Solve the proximal subproblem of second-order cone} \label{append-prox-soc}
	We shall solve (\ref{cp-prox}) in two cases of second-order cone $\mathcal{K}$: SOC and RSOC. The subproblem (\ref{cp-prox}) is
\begin{equation}
	\bx^{*} = \argmin_{\bx\in\text{int}\mathcal{K}} \lambda F(\bx) + \frac{1}{2}\|\bx - \bzeta\|_{2}^{2}
\end{equation}
where $\lambda > 0$, $\zeta$ is a column vector and $F(\bx)$ is the barrier function of $\mathcal{K}$. The optimality condition is
\begin{equation}
    \lambda \nabla F(\bx) + \bx - \bzeta = \bzero \label{prox-oc}
\end{equation}

\subsection{Regular second-order cone}
Let $\bx = (t,\bar{\bx})\in\mathbb{R}\times\mathbb{R}^{n_{\bar{\bx}}}$ and $\mathcal{K} = \left\lbrace \bx: t \geq \sqrt{\bar{\bx}^{T}\bar{\bx}}\right\rbrace$. Notice that $F(\bx) = - \log \left( t^{2} - \bar{\bx}^{T}\bar{\bx} \right)$.

In the view of the optimality condition \eqref{prox-oc}, we need to solve
\begin{equation*}
	- \frac{\lambda}{t^{2} - \bar{\bx}^{T}\bar{\bx}}
	\begin{bmatrix}
		2t \\
		-2\bar{\bx}
	\end{bmatrix} + \begin{bmatrix} t \\ \bar{\bx} \end{bmatrix} - \begin{bmatrix} \zeta_{t} \\ \bzeta_{\bar{\bx}} \end{bmatrix} = \bzero.
\end{equation*}
For brevity, let us denote
\begin{equation}
	\delta = t^{2} - \bar{\bx}^{T}\bar{\bx} > 0 \label{soc-delta}
\end{equation}
we have
\begin{align}
	(-2\lambda + \delta) t &= \delta\zeta_{t}, \label{soc-1}\\
	(2\lambda + \delta) \bar{\bx} &= \delta\bzeta_{\bar{\bx}}. \label{soc-2}
\end{align}
In view of the relation \eqref{soc-1},  we consider two cases.

1) When $\delta=2\lambda$, from (\ref{soc-1}) and (\ref{soc-2}) we have $\zeta_{t} = 0, \bar{\bx}^* = \frac{1}{2}\bzeta_{\bar{\bx}}$. Plugging the value of $\bar{\bx}^*$ in \eqref{soc-delta} and combining it with $\delta=2\lambda$, we have
\begin{equation*}
	2\lambda = t^{2} - \frac{1}{4}\bzeta_{\bar{\bx}}^{T}\bzeta_{\bar{\bx}}\ \Rightarrow\ t^* = \sqrt{2\lambda + \frac{1}{4}\bzeta_{\bar{\bx}}^{T}\bzeta_{\bar{\bx}}}.
\end{equation*}

2) When $\delta\neq 2\lambda$, from (\ref{soc-1}) and (\ref{soc-2}) we have
\begin{equation}
	\begin{cases}
		t = \frac{\delta}{\delta - 2\lambda}\zeta_{t} \\
		\bar{\bx} = \frac{\delta}{\delta + 2\lambda}\bzeta_{\bar{\bx}}
	\end{cases}. \label{soc-t-x}
\end{equation}
Substitute $t, \bar{\bx}$ in \eqref{soc-delta}, we have the equation of $\delta$
\begin{align*}
	&\delta = \frac{\delta^{2}}{(\delta-2\lambda)^{2}}\zeta_{t}^{2} - \frac{\delta^{2}}{(\delta+2\lambda)^{2}}\bzeta_{\bar{\bx}}^{T}\bzeta_{\bar{\bx}} \\
	\Rightarrow\ & (\delta+2\lambda)^{2}(\delta-2\lambda)^{2} = \delta (\delta+2\lambda)^{2} \zeta_{t}^{2}  - \delta (\delta-2\lambda)^{2}\bzeta_{\bar{\bx}}^{T}\bzeta_{\bar{\bx}}.
\end{align*}
Let $\rho=\frac{\delta}{\lambda}\in(0,2)\cup(2,+\infty)$, then we have the equation of $\rho$
\begin{align*}
	& (\rho + 2)^{2}(\rho - 2)^{2} = \frac{\rho}{\lambda}(\rho+2)^{2}\zeta_{t}^{2} - \frac{\rho}{\lambda}(\rho-2)^{2}\bzeta_{\bar{\bx}}^{T}\bzeta_{\bar{\bx}} \\
	\Leftrightarrow\ & \rho^{4} - 8\rho^{2} + 16 = \frac{\zeta_{t}^{2} - \bzeta_{\bar{\bx}}^{T}\bzeta_{\bar{\bx}}}{\lambda}\rho^{3} + \frac{4\zeta_{t}^{2}+ 4\bzeta_{\bar{\bx}}^{T}\bzeta_{\bar{\bx}}}{\lambda}\rho^{2} + \frac{4\zeta_{t}^{2} - 4\bzeta_{\bar{\bx}}^{T}\bzeta_{\bar{\bx}}}{\lambda}\rho \\
	\Leftrightarrow\ & \left(\rho^{2} + \frac{16}{\rho^{2}} \right) - \frac{\zeta_{t}^{2} - \bzeta_{\bar{\bx}}^{T}\bzeta_{\bar{\bx}}}{\lambda}\left(\rho + \frac{4}{\rho} \right) - \frac{4\zeta_{t}^{2}+ 4\bzeta_{\bar{\bx}}^{T}\bzeta_{\bar{\bx}}}{\lambda} - 8 = 0 \\
	\Leftrightarrow\ & \left(\rho + \frac{4}{\rho} \right)^{2} - \frac{\zeta_{t}^{2} - \bzeta_{\bar{\bx}}^{T}\bzeta_{\bar{\bx}}}{\lambda}\left(\rho + \frac{4}{\rho} \right) - \frac{4\zeta_{t}^{2}+ 4\bzeta_{\bar{\bx}}^{T}\bzeta_{\bar{\bx}}}{\lambda} - 16 = 0.
\end{align*}
Let $\gamma = \rho + \frac{4}{\rho} > 4$, then we have the equation of $\gamma$
\begin{equation}
	\gamma^{2} - \frac{\zeta_{t}^{2} - \bzeta_{\bar{\bx}}^{T}\bzeta_{\bar{\bx}}}{\lambda}\gamma - \frac{4\zeta_{t}^{2}+ 4\bzeta_{\bar{\bx}}^{T}\bzeta_{\bar{\bx}}}{\lambda} - 16 = 0. \label{soc-eq-gamma}
\end{equation}
It is easy to prove that the equation above must have and only have one solution no less than 4, so we choose the larger solution
\begin{equation*}
	\gamma = \frac{\frac{\zeta_{t}^{2} - \bzeta_{\bar{\bx}}^{T}\bzeta_{\bar{\bx}}}{\lambda} + \sqrt{\left( \frac{\zeta_{t}^{2} - \bzeta_{\bar{\bx}}^{T}\bzeta_{\bar{\bx}}}{\lambda}\right)^{2} + 4\left( \frac{4\zeta_{t}^{2}+ 4\bzeta_{\bar{\bx}}^{T}\bzeta_{\bar{\bx}}}{\lambda} + 16 \right) }}{2}
\end{equation*}
and
\begin{equation*}
	\rho_{1} = \frac{\gamma - \sqrt{\gamma^{2}-16}}{2} < 2, \quad \rho_{2} = \frac{\gamma + \sqrt{\gamma^{2}-16}}{2} > 2.
\end{equation*}
In order to make $t \geq 0$, from (\ref{soc-t-x}), $\rho$ should satisfy
\begin{equation*}
	(\rho-2)\zeta_{t} \geq 0.
\end{equation*}
Thus if $\zeta_{t} > 0$, we choose $\rho_{2}$; if $\zeta_{t} < 0$, we choose $\rho_{1}$. When $\zeta_{t}\neq 0$, from (\ref{soc-t-x}), we have
\begin{equation*}
	t^* = \frac{\rho}{\rho - 2}\zeta_{t}, \quad
	\bar{\bx}^* = \frac{\rho}{\rho+2}\bzeta_{\bar{\bx}}.
\end{equation*}

\subsection{Rotated second-order cone}
Let $\bx = (\eta,\nu,\bar{\bx}) \in \mathbb{R}_{+}\times\mathbb{R}_{+}\times\mathbb{R}^{n_{\bar{\bx}}}$ and $\mathcal{K} = \left\lbrace \bx: \eta\nu \geq \frac{1}{2}\bar{\bx}^{T}\bar{\bx}, \eta\geq 0,\nu \geq 0 \right\rbrace$. Notice that $F(\bx) = -\log \left(\eta\nu - \frac{1}{2}\bar{\bx}^{T}\bar{\bx} \right)$.

In view of the optimality condition \eqref{prox-oc}, we have
\begin{equation*}
	- \frac{\lambda}{\eta\nu - \frac{1}{2}\bar{\bx}^{T}\bar{\bx}}
	\begin{bmatrix}
		\nu \\
		\eta \\
		-\bar{\bx}
	\end{bmatrix} + \begin{bmatrix} \eta \\ \nu \\ \bar{\bx} \end{bmatrix} - \begin{bmatrix} \zeta_{\eta} \\ \zeta_{\nu} \\ \bzeta_{\bar{\bx}} \end{bmatrix} = \bzero.
\end{equation*}
For brevity, let us denote 
\begin{equation}
	\delta = \eta\nu - \frac{1}{2}\bar{\bx}^{T}\bar{\bx} > 0, \label{rsoc-delta}
\end{equation}
we have
\begin{align}
	-\lambda\nu + \delta \eta &= \delta\zeta_{\eta}, \label{rsoc-1}\\
	-\lambda\eta + \delta \nu &= \delta\zeta_{\nu}, \label{rsoc-2}\\
	(\lambda + \delta) \bar{\bx} &= \delta\bzeta_{\bar{\bx}}. \label{rsoc-3}
\end{align}
Due to the relation \eqref{rsoc-1} and \eqref{rsoc-2} we consider two cases.

1) When $\delta = \lambda$, from (\ref{rsoc-1}) and (\ref{rsoc-2}) we have $\zeta_{\eta} + \zeta_{\nu} = 0$. From (\ref{rsoc-3}) we have
$\bar{\bx}^* = \frac{1}{2}\bzeta_{\bar{\bx}}$. 
Plugging this value of $\bar{\bx}^*$ in \eqref{rsoc-delta} and combining it with $\delta=\lambda$ and \eqref{rsoc-1} gives
\begin{equation*}
	\begin{cases}
		\eta - \nu = \zeta_{\eta} \\
		\lambda = \eta\nu - \frac{1}{8}\bzeta_{\bar{\bx}}^{T}\bzeta_{\bar{\bx}}
	\end{cases}\ \Rightarrow\ 
	\begin{cases}
		\eta^* = \frac{\zeta_{\eta} + \sqrt{\zeta_{\eta}^{2} + 4\left(\lambda + \frac{1}{8}\bzeta_{\bar{\bx}}^{T}\bzeta_{\bar{\bx}} \right) }}{2} \\
		\nu^* = \frac{-\zeta_{\eta} + \sqrt{\zeta_{\eta}^{2} + 4\left(\lambda + \frac{1}{8}\bzeta_{\bar{\bx}}^{T}\bzeta_{\bar{\bx}} \right) }}{2}
	\end{cases}.
\end{equation*}

2) When $\delta\neq\lambda$, from (\ref{rsoc-1}) and (\ref{rsoc-2}) we have
\begin{equation}
	\begin{cases}
		\eta + \nu = \frac{\delta}{\delta - \lambda}(\zeta_{\eta} + \zeta_{\nu}) \\
		\eta - \nu = \frac{\delta}{\delta + \lambda}(\zeta_{\eta} - \zeta_{\nu})
	\end{cases} \ \Rightarrow\ 
	\begin{cases}
		\eta = \frac{\delta^{2}\zeta_{\eta} + \delta\lambda\zeta_{\nu}}{(\delta+\lambda)(\delta-\lambda)} \\
		\nu = \frac{\delta\lambda\zeta_{\eta} + \delta^{2}\zeta_{\nu}}{(\delta+\lambda)(\delta-\lambda)}
	\end{cases}. \label{rsoc-eta-nu}
\end{equation}
From (\ref{rsoc-3}) we have $\bar{\bx} = \frac{\delta\bzeta_{\bar{\bx}}}{\delta+\lambda}$. Substitute $\eta, \nu, \bar{\bx}$ in (\ref{rsoc-delta}), we have the equation of $\delta$
\begin{align*}
	&\delta = \frac{\delta^{2}}{(\delta+\lambda)^{2}(\delta-\lambda)^{2}}(\delta\zeta_{\eta} + \lambda\zeta_{\nu})(\lambda\zeta_{\eta} + \delta\zeta_{\nu}) - \frac{1}{2}\frac{\delta^{2}}{(\delta+\lambda)^{2}}\bzeta_{\bar{\bx}}^{T}\bzeta_{\bar{\bx}} \\
	\Rightarrow\ & 2(\delta+\lambda)^{2}(\delta-\lambda)^{2} = 2\delta(\delta\zeta_{\eta} + \lambda\zeta_{\nu})(\lambda\zeta_{\eta} + \delta\zeta_{\nu}) - \delta (\delta-\lambda)^{2}\bzeta_{\bar{\bx}}^{T}\bzeta_{\bar{\bx}}.
\end{align*}

Let $\rho=\frac{\delta}{\lambda}\in (0,1)\cup (1,+\infty)$, then we have the equation of $\rho$
\begin{align*}
	& 2(\rho+1)^{2}(\rho-1)^{2} = 2\frac{\rho}{\lambda}(\zeta_{\eta}\rho+\zeta_{\nu})(\zeta_{\eta}+\zeta_{\nu}\rho) - \frac{\rho}{\lambda}(\rho-1)^{2}\bzeta_{\bar{\bx}}^{T}\bzeta_{\bar{\bx}} \\
	\Leftrightarrow\ &2(\rho^{4} - 2\rho^{2} + 1) = \frac{2\zeta_{\eta}\zeta_{\nu} - \bzeta_{\bar{\bx}}^{T}\bzeta_{\bar{\bx}}}{\lambda} \rho^{3} + \frac{2\zeta_{\eta}^{2} + 2\zeta_{\nu}^{2} + 2\bzeta_{\bar{\bx}}^{T}\bzeta_{\bar{\bx}}}{\lambda} \rho^{2} + \frac{2\zeta_{\eta}\zeta_{\nu} - \bzeta_{\bar{\bx}}^{T}\bzeta_{\bar{\bx}}}{\lambda} \rho \\
	\Leftrightarrow\ &2\left( \rho^{2} + \frac{1}{\rho^{2}} \right) - \frac{2\zeta_{\eta}\zeta_{\nu} - \bzeta_{\bar{\bx}}^{T}\bzeta_{\bar{\bx}}}{\lambda}\left( \rho + \frac{1}{\rho} \right) - \frac{2\zeta_{\eta}^{2} + 2\zeta_{\nu}^{2} + 2\bzeta_{\bar{\bx}}^{T}\bzeta_{\bar{\bx}}}{\lambda} - 4 = 0 \\
	\Leftrightarrow\ &2\left( \rho + \frac{1}{\rho} \right)^{2} - \frac{2\zeta_{\eta}\zeta_{\nu} - \bzeta_{\bar{\bx}}^{T}\bzeta_{\bar{\bx}}}{\lambda}\left( \rho + \frac{1}{\rho} \right) - \frac{2\zeta_{\eta}^{2} + 2\zeta_{\nu}^{2} + 2\bzeta_{\bar{\bx}}^{T}\bzeta_{\bar{\bx}}}{\lambda} - 8 = 0.
\end{align*}

Let $\gamma = \rho + \frac{1}{\rho} > 2$, then we have the equation of $\gamma$
\begin{equation}
	\gamma^{2} - \frac{2\zeta_{\eta}\zeta_{\nu} - \bzeta_{\bar{\bx}}^{T}\bzeta_{\bar{\bx}}}{2\lambda}\gamma - \frac{\zeta_{\eta}^{2} + \zeta_{\nu}^{2} + \bzeta_{\bar{\bx}}^{T}\bzeta_{\bar{\bx}}}{\lambda} - 4 = 0. \label{rsoc-eq-gamma}
\end{equation}
It is easy to prove that the equation above must have and only have one solution no less than 2, so we choose the larger solution
\begin{equation*}
	\gamma = \frac{\frac{2\zeta_{\eta}\zeta_{\nu} - \bzeta_{\bar{\bx}}^{T}\bzeta_{\bar{\bx}}}{2\lambda} + \sqrt{\left( \frac{2\zeta_{\eta}\zeta_{\nu} - \bzeta_{\bar{\bx}}^{T}\bzeta_{\bar{\bx}}}{2\lambda}\right)^{2} + 4\left( \frac{\zeta_{\eta}^{2} + \zeta_{\nu}^{2} + \bzeta_{\bar{\bx}}^{T}\bzeta_{\bar{\bx}}}{\lambda} + 4 \right) }}{2}
\end{equation*}
and then we can get the possible values of $\rho$
\begin{equation*}
	\rho_{1} = \frac{\gamma - \sqrt{\gamma^{2}-4}}{2} < 1, \quad \rho_{2} = \frac{\gamma + \sqrt{\gamma^{2}-4}}{2} > 1.
\end{equation*}
In order to ensure $\eta\geq0$ and $\nu \geq 0$, from \eqref{rsoc-eta-nu} $\rho$ should satisfy
\[
	(\rho-1)(\rho\zeta_{\eta} + \zeta_{\nu}) \geq 0, \ \text{ and }
	(\rho-1)(\zeta_{\eta} + \rho\zeta_{\nu}) \geq 0.
\]
Therefore, if $\zeta_{\eta} + \zeta_{\nu} > 0$, we set $\rho=\rho_{2}$; if $\zeta_{\eta} + \zeta_{\nu} < 0$, we set $\rho=\rho_{1}$. When $\zeta_{\eta}+\zeta_{\nu}\neq 0$, from \eqref{rsoc-eta-nu} we have
\begin{equation*}
		\eta^* = \frac{\rho^{2}\zeta_{\eta} + \rho\zeta_{\nu}}{(\rho+1)(\rho-1)}, \quad
		\nu^* = \frac{\rho\zeta_{\eta} + \rho^{2}\zeta_{\nu}}{(\rho+1)(\rho-1)},\quad
		\bar{\bx}^* = \frac{\rho\bzeta_{\bar{\bx}}}{\rho+1}.
\end{equation*}

\section{Specialized linear system solvers and two equivalent formulations for LASSO and SVM}
\label{sec:linsys}

To further improve the performance of ABIP on several important applications, we provide specialized linear system solvers that can effectively exploit the problem structure. 

\subsection{LASSO}\label{subsec:lasso}
Consider the LASSO problem:
	\begin{equation*}
		\min_{\bx}\ \|\tilde{\bA}\bx-\bb\|_{2}^{2} + \lambda\|\bx\|_{1}
	\end{equation*}
	where $\tilde{\bA} \in \mathbb{R}^{m\times n}$. 
 
 \textbf{SOCP} reformulation for ABIP is:
	\begin{align*}
		\min\ & \begin{bmatrix}
			0 \\ 2 \\ \bzero_{m\times 1} \\ \lambda\cdot\mathbf{1}_{2n\times 1}
		\end{bmatrix}^{T}\begin{bmatrix}
			w \\ z \\ \by \\ \bx^{+} \\ \bx^{-}
		\end{bmatrix} \\
		\text{s.t.}\ &\begin{bmatrix}
			1 & 0 & \bzero_{1\times (m+2n)} \\
			\bzero_{m\times 1} & \bzero_{m\times 1} & \bI_{m} & \tilde{\bA} & -\tilde{\bA}
		\end{bmatrix} \begin{bmatrix}
			w \\ z \\ \by \\ \bx^{+} \\ \bx^{-}
		\end{bmatrix} = \begin{bmatrix}
			1 \\ \bb
		\end{bmatrix} \\
		& (w,z,\by)^T \in \text{RSOC}_{2+m} \\
		& \bx^{+},\bx^{-}\geq \bzero
	\end{align*}
	As discussed above, our aim is to solve $(\bI + \bA\bA^{T}) \bx = \bb$, and for LASSO, 
	\begin{equation*}
		\bA\bA^{T} = \begin{bmatrix}
			1 & \bzero_{1\times m} \\
			\bzero_{m\times 1} & \bI_{m}+2\tilde{\bA}\tilde{\bA}^T
		\end{bmatrix},\quad
	\end{equation*}
	so we need to solve $(\bI_{m} + 2\tilde{\bA}\tilde{\bA}^T) \bx = \bb$, then we may reduce the dimension of factorization.
	\begin{enumerate}
	    \item If $n > m$, we directly perform cholesky or LDL factorization to $\bI_{m} + 2\tilde{\bA}\tilde{\bA}^T$. In this case, we reduce the dimension of factorization from $2m+2n+3$ to $m$.
	    \item If $n\leq m$, we apply the Sherman-Morrison-Woodbury formula again to $\bI_{m} + 2\tilde{\bA}\tilde{\bA}^T$, we can get:
	    \begin{align*}
	    (\bI_{m} + 2\tilde{\bA}\tilde{\bA}^T)^{-1} = & \bI_{m} - 2\tilde{\bA}(\bI_{m}+2\tilde{\bA}^T\tilde{\bA})^{-1}\tilde{\bA}^{T}\\
	    = & \bI_{m}-\tilde{\bA}(0.5\bI_{m}+\tilde{\bA}^T\tilde{\bA})^{-1}\tilde{\bA}^{T}
	    \end{align*}
	    Then, we only need to perform Cholesky or LDL factorization to $0.5\bI_{m}+\tilde{\bA}^T\tilde{\bA}$. In this case, we reduce the dimension of factorization from $2m+2n+3$ to $n$.
	\end{enumerate}
	In a word, we only need to perform matrix factorization in the dimension of $O(\min\{m,n\})$ instead of $O(m+n)$.

 \textbf{QP} reformulation for ABIP is:
$$\min \ \hat{\bx}^T\hat{\bA}^T\hat{\bA}\hat{\bx} + (\lambda\cdot \mathbf{1}_{2n}-2\hat{\bA}^T\bb)^T\hat{\bx}$$
    where $\hat{\bx}=\begin{bmatrix}
        \bx^+\\
        \bx^-
    \end{bmatrix} \geq \bzero,\ \bx = \bx^+ - \bx^- ,\ \hat{\bA} = [\tilde{\bA}\quad -\tilde{\bA}] $

\subsection{Support vector machines}\label{subsec:svm}
Consider the Support Vector Machine problem with training data $(\bx_i,y_i)_{1\le i\le n}$:
\begin{equation}\label{eq:svm-prob}
	\begin{aligned}
	\min \ & \frac{1}{2}\|\tilde{\bw}\|^{2}+C \sum_{i=1}^{m} \xi_{i} \\
	\text { s.t. } & y_{i} \bx_{i}^{T} \tilde{\bw}+y_{i} \tilde{b}+\xi_{i} \geqslant 1 \quad i=1, \ldots, m \\
	& \xi_{i} \geqslant 0 \quad i=1, \ldots, m
	\end{aligned}
\end{equation}	
where each $\bx_i \in \mathbb{R}^n$. 
For brevity, let $\bX=[\bx_1,\bx_2,\ldots,\bx_m]^T\in\mathbb{R}^{m\times n} $ and
 $\by=[y_1,y_2,\ldots,y_n]^T\in\mathbb{R}^n$ represent the feature matrix and the label vector, respectively.

\textbf{SOCP} reformulation of SVM~\eqref{eq:svm-prob} for ABIP is:
\begin{align*}
	\min &\quad \nu+C \sum_{i=1}^{l} \xi_{i} \\
	\text { s.t. } & \quad \eta=1 \\
	& \quad \tilde{\bA} \tilde{\bw}^{+}-\tilde{\bA} \tilde{\bw}^{-}+\tilde{b}^{+} \by-\tilde{b}^{-} \by+\xi-\bt=\mathbf{e}_m \\
	& \quad \tilde{\bw}=\tilde{\bw}^{+}-\tilde{\bw}^{-} \\
	& \quad \xi, \tilde{\bw}^{+}, \tilde{\bw}^{-}, \tilde{b}^{+}, \tilde{b}^{-}, \bt \geqslant 0 \\
	& \quad \eta \nu-\frac{1}{2}\|\tilde{\bw}\|^{2} \ge 0,
\end{align*}
where $\mathbf{e}_m$ is the vector that all elements are one, and $\tilde{\bA} = \textup{diag}(\by)\cdot \bX$. Our formulation introduces new variables $\tilde{\bw}^+$ and $\tilde{\bw}^-$, which puts the matrix $\tilde{\bA}$ into the columns of positive orthant. An empirical advantage of such formulation is that ABIP can scale each column independently. Note that, however, such scaling can not be accomplished for the initial formulation, as the matrix $\tilde{\bA}$ is in the columns of the RSOC and can only be scaled as a whole.\\
Consequently, we reformulate SVM to align with the ABIP standard input
\begin{equation*}\label{prob:conic}
    \begin{aligned}
        \min & \quad \bc^{T}\bx \\
    \quad	\text{s.t.} &\quad \bA\bx = \bb \\
        &\quad \bx \in \mathcal{K}
    \end{aligned}
\end{equation*}
where
\begin{equation*}
	\bA = \begin{bmatrix}
		1 &  \\
		& & & \tilde{\bA} & \by & -\tilde{\bA} & -\by & \bI_m & -\bI_m \\
		& & \bI_n & -\bI_n &  & \bI_n 
	\end{bmatrix},
\end{equation*}
$$
\bx=\begin{bmatrix}
\eta \\
\nu \\
\tilde{\bw} \\
\tilde{\bw}^{+} \\
\tilde{b}^{+} \\
\tilde{\bw}^{-} \\
\tilde{b}^{-} \\
\bxi \\
\bt
\end{bmatrix}, 
\bc=\begin{bmatrix}
0 \\
1 \\
\bzero_{n} \\
\bzero_{n} \\
0 \\
\bzero_{n} \\
0 \\
C \mathbf{e}_m \\
\bzero_{m}
\end{bmatrix},
\ \text{and}\ 
\bb = \begin{bmatrix}
    1 \\ \mathbf{e}_{m} \\ \bzero_{n}
\end{bmatrix},
\begin{bmatrix}
    \eta \\
    \nu \\
    \tilde{\bw}
\end{bmatrix}
\in \text{RSOC}_{2+n},
\begin{bmatrix}
    \tilde{\bw}^{+} \\
    \tilde{b}^{+} \\
    \tilde{\bw}^{-} \\
    \tilde{b}^{-} \\
    \bxi \\
    t
\end{bmatrix}
\in \mathbb{R}^{2n+2m+2}_+.
$$
Therefore, we need to perform the LDL factorization of the following matrix:
$$
\bK=\begin{bmatrix}
\bI_{m+n+1} & \bA \\
\bA^{T} & -\bI_{2 m+3 n+4}
\end{bmatrix},
$$
where
\begin{equation*}
	\bA = \begin{bmatrix}
		1 &  \\
		& & & \tilde{\bA} & \by & -\tilde{\bA} & -\by & \bI_m & -\bI_m \\
		& & \bI_n & -\bI_n &  & \bI_n.
	\end{bmatrix}.
\end{equation*}
Note that there is an obvious factorization form of $\bK$:
$$
\bK = 
\left[\begin{array}{cc}
\bI_{m+n+1} & -\bA \\
& \bI_{2 m+3 n+4}
\end{array}\right]\left[\begin{array}{cc}
\bI_{m+n+1}+\bA \bA^{T} & \\
& -\bI_{2 m+3 n+4}
\end{array}\right]\left[\begin{array}{cc}
\bI_{m+n+1} & \\
-\bA^{T} & \bI_{2 m+3 n+4}
\end{array}\right]
$$
Then, it suffices to factorize the matrix:
\begin{align*}
& \bI_{m+n+1}+\bA \bA^{T} 
\\ &= 
\begin{bmatrix}
2 & & \\
& 2 \tilde{\bA} \tilde{\bA}^{T}+2 \by \by^{T}+3 \bI_{m} & -2 \tilde{\bA}\\
& -2 \tilde{\bA}^{T} & 4 \bI_{n}
\end{bmatrix}
\\& = 
\begin{bmatrix}
	2 & & \\
	& 2 \tilde{\bA} \tilde{\bA}^{T}+2 \by \by^{T}+\bF & \bX \\
	& \bX^{T} & \bG
\end{bmatrix}
\\& = 
\begin{bmatrix}
	1 & & \\
	& \bI_{m} & \bX \bG^{-1} \\
	& & \bI_{n}
\end{bmatrix}
\begin{bmatrix}
	2 & & & \\
	& 2 \tilde{\bA} \tilde{\bA}^{T}+2 \by \by^{T}+\bF-\bX \bG^{-1} \bX^{T} & \\
	& & \bG
\end{bmatrix}
\begin{bmatrix}
	1 & & \\
	& \bI_{m} & \\
& \bG^{-1} \bX^{T} & \bI_{n}
\end{bmatrix}
\\& = 
\begin{bmatrix}
	1 & & \\
	& \bI_{m} & -2 \tilde{\bA} \bE \bG^{-1} \\
	& & \bI_{n}
\end{bmatrix}
\begin{bmatrix}
	2 & & \\
	& 2 \tilde{\bA} \tilde{\bA}^{T}+2 \by \by^{T}+\bF-4 \tilde{\bA} \bE \bG^{-1} \bE \tilde{\bA}^{T} & \\
	& & \bG
\end{bmatrix}
\begin{bmatrix}
	1 & & \\
	& \bI_{m} & \\
	& -2 \bG^{-1} \bE \tilde{\bA}^{T} & \bI_{n}
\end{bmatrix}
\end{align*}
where $\bX=-2 \tilde{\bA}, \bF=3 \bI_{m}, \bG=4 \bI_{n}$. Let 
\begin{equation*}
	\bH = 
	\begin{bmatrix}
		2 \bI_{n}-4 \bG^{-1} & \\ & 2
	\end{bmatrix}
	\quad \text{and}
	\quad
	\bM = 
	\begin{bmatrix}
		\tilde{\bA} & \by
	\end{bmatrix}.
\end{equation*}
It follows that  $\bH$ is a diagonal matrix and
$$2 \tilde{\bA} \tilde{\bA}^{T}+2 \by \by^{T}+\bF-4 \tilde{\bA} \bG^{-1} \tilde{\bA}^{T}=\bF+\bM \bH \bM^{T}.$$
We consider two cases for factorization. 1) When $m \leq n+1$, we directly factorize $\bF+\bM \bH \bM^{T}$. 2) When $m > n+1$, we apply the Sherman-Morrison-Woodbury formula to obtain:
$$
\left(\bF+\bM \bH \bM^{T}\right)^{-1}=\bF^{-1}-\bF^{-1} \bM\left(\bH^{-1}+\bM^{T} \bF^{-1} \bM\right)^{-1} \bM^{T} \bF^{-1}.
$$
In the above case, we factorize $\bH^{-1}+\bM^{T} \bF^{-1} \bM$. Therefore, we reduce the dimension of 
matrix factorization from $3m+4n+5$ to $\min\{m,n+1\}$.

\textbf{QP} reformulation of SVM~\eqref{eq:svm-prob} for ABIP is:
\begin{align*}
    \min & \quad \frac{\lambda}{2} \tilde{\bx}^T\begin{bmatrix}
        I_{n} &  \\
        & \mathbf{0}_{2m+1}
    \end{bmatrix} \tilde{\bx} + 
    \frac{1}{m} \begin{bmatrix}
        &\textbf{0}_{n+1}\\
        &\mathbf{1}_{m}\\
        &\textbf{0}_{m}
    \end{bmatrix}^T \tilde{\bx}\\
    \text{s.t.}& \quad \begin{bmatrix}
        \textrm{diag}(\by)\cdot \bX&\by&\bI_{m}&-\bI_{m}
    \end{bmatrix}\tilde{\bx} = \mathbf{1}_{m}\\
    &\quad \tilde{\bx} = \begin{bmatrix}
        \bw\\\bb\\\xi\\s
    \end{bmatrix},\quad \xi,s \geq 0.
\end{align*}
\end{document}